%% file: z.tex
\newif\ifsiopt

\sioptfalse 

\RequirePackage{fix-cm}
\ifsiopt
        \documentclass[review,nohypdvips]{siamart1116}
\else
        \documentclass[smallextended]{svjour3}
        \smartqed
\fi

\usepackage{graphicx}
\usepackage{amsmath,amssymb}
\usepackage{array}
\usepackage{enumitem}
\usepackage{xcolor}
\usepackage{bm}
\usepackage{multirow}
\usepackage{xspace}
\usepackage{pbox}
\usepackage{hyphenat}
\input def,command

\newcommand{\FO} {FSFOM\xspace}
\newcommand{\m} {\xmath{m}\xspace}

\renewcommand{\Theta} {\xmath{\Omega}} 

\newcommand{\vTheta} {\bmath{{\Theta}}}
\newcommand{\vtheta} {\bmath{{\theta}}}

\newcommand{\vT} {\bmath{{T}}}
\newcommand{\vt} {\bmath{{t}}}

\newcommand{\bmlam} {\bmath{{\lambda}}}
\newcommand{\bmtau} {\bmath{{\tau}}}

\newcommand{\bmh} {\bmath{h}}
\newcommand{\g} {\bmath{g}}
\newcommand{\G} {\bmath{G}}
\renewcommand{\del} {\bmath{\delta}}
\renewcommand{\u} {\bmath{u}}
\newcommand{\e} {\bmath{e}}
\newcommand{\x} {\bmath{x}}
\newcommand{\y} {\bmath{y}}
\newcommand{\z} {\bmath{z}}
\newcommand{\bmnu} {\bmath{\nu}}
\renewcommand{\S} {\bmath{S}}
\newcommand{\A} {\bmath{A}}
\newcommand{\B} {\bmath{B}}
\newcommand{\C} {\bmath{C}}
\newcommand{\D} {\bmath{D}}
\newcommand{\bmbeta} {\bmath{{\beta}}}

\newcommand{\Zero} {\bmath{0}}
\newcommand{\Q} {\bmath{Q}}
\newcommand{\p} {\bmath{p}}
\newcommand{\rr} {\bmath{r}}

\newcommand{\cC} {\xmath{\mathcal{C}_L^{1,1}(\Reals^d)}}
\newcommand{\cF} {\xmath{\mathcal{F}_L(\Reals^d)}}
\newcommand{\cQ} {\xmath{\mathcal{Q}_L(\Reals^d)}}
\newcommand{\cL} {\xmath{\mathcal{L}}}
\newcommand{\Reals} {\xmath{\mathbb{R}}}

\newcommand{\hki} {\xmath{h_{i,k}}}
\newcommand{\hkip} {\xmath{h_{i+1,k}}}
\newcommand{\hkj} {\xmath{h_{j,k}}}
\newcommand{\hiip} {\xmath{h_{i+1,i}}}
\newcommand{\himi} {\xmath{h_{i,i-1}}}
\newcommand{\himip} {\xmath{h_{i+1,i-1}}}

\newcommand{\hkn} {\xmath{h_{n,k}}}
\newcommand{\hknp} {\xmath{h_{n+1,k}}}
\newcommand{\hkl} {\xmath{h_{l,k}}}

\newcommand{\fNh} {\xmath{\floor{\Frac{N}{2}}}}

\newcommand{\fNt} {\xmath{\floor{\Frac{2N}{3}}}}

\numberwithin{equation}{section}
\ifsiopt
	\numberwithin{theorem}{section}
\else
	\usepackage{fullpage}
\fi
\usepackage{hyperref}
\hypersetup{
        linktoc=page,
        colorlinks=true,
        linkcolor=blue,
        citecolor=red,
        filecolor=magenta,
        urlcolor=cyan
}

\newcommand{\TheTitle}{Generalizing the optimized gradient method
	for smooth convex minimization}
\newcommand{\TheAuthors}{D. Kim, and J. A. Fessler}

\ifsiopt

        \headers{Generalizing the optimized gradient method}{\TheAuthors}

        \title{{\TheTitle}\thanks{Submitted to the editors \today.
        \funding{This research was supported in part by NIH grant U01 EB018753.}}}

        \author{
        Donghwan Kim\thanks{Dept. of Electrical Engineering and Computer Science,
                University of Michigan, Ann Arbor, MI 48109 USA
                (\email{kimdongh@umich.edu}, \email{fessler@umich.edu})}
        \and
        Jeffrey A. Fessler\footnotemark[2]
        }

        \ifpdf
        \hypersetup{
        pdftitle={\TheTitle}
        pdfauthor={\TheAuthors}
        }
        \fi

\else
        \title{\TheTitle
                \thanks{This research was supported in part by NIH grant U01 EB018753.}}
        \author{Donghwan Kim \and Jeffrey A. Fessler}
        \institute{Donghwan Kim \and Jeffrey A. Fessler \at
              Dept. of Electrical Engineering and Computer Science,
                University of Michigan, Ann Arbor, MI 48109 USA \\
              \email{kimdongh@umich.edu, fessler@umich.edu}           
        }
        \date{Date of current version: \today} 
\fi

\begin{document}
\interfootnotelinepenalty=10000

\maketitle

\begin{abstract}
This paper generalizes 
the optimized gradient method (OGM)~\cite{drori:14:pof,kim:16:ofo,kim:17:otc}
that achieves
the optimal worst-case cost function bound of first-order methods
for smooth convex minimization~\cite{drori:17:tei}.
Specifically,
this paper studies a generalized formulation of OGM
and analyzes its worst-case rates
in terms of both the function value and the norm of the function gradient.
This paper also develops a new algorithm called OGM-OG
that is in the generalized family of OGM and
that has the best known analytical worst-case bound
with rate $O(1/N^{1.5})$
on the decrease of the gradient norm
among fixed-step first-order methods.
This paper also proves
that Nesterov's fast gradient method~\cite{nesterov:83:amf,nesterov:05:smo}
has an $O(1/N^{1.5})$ worst-case gradient norm rate
but with constant larger than OGM-OG.
The proof is based on the worst-case analysis 
called Performance Estimation Problem in~\cite{drori:14:pof}.
%
\end{abstract}

\section{Introduction}
\label{intro}

First-order methods are favorable
for solving large-scale problems
because their computational complexity per iteration
depends mildly on the problem dimension.
In particular, Nesterov's fast gradient method (FGM)
\cite{nesterov:83:amf,nesterov:05:smo} 
achieves the optimal worst-case rate $O(1/N^2)$
for decreasing smooth convex functions after $N$ iterations
\cite{nesterov:04},
and thus has been widely used in (large-scale) applications.
Recently,
the optimized gradient method (OGM)~\cite{drori:14:pof,kim:16:ofo,kim:17:otc}
was found to achieve the optimal worst-case cost function bound
of first-order methods 
(with either fixed-step or adaptive-step approaches)
for smooth convex minimization in~\cite{drori:17:tei},
whereas FGM achieves that bound only up to constant.\footnote{
There is a backtracking line-search version of FGM~\cite{nesterov:04}
that also achieves the optimal worst-case function bound up to constant,
which is sometimes more useful than the fixed-step FGM in practice.
However, such backtracking line-search version of OGM 
with a fast worst-case bound
is yet unknown (unlike the fixed-step OGM~\cite{drori:14:pof,kim:16:ofo,kim:17:otc}),
while recently an exact line-search version of OGM 
is developed in~\cite{drori:18:efo}.
}

Building upon~\cite{drori:14:pof,kim:16:ofo,kim:17:otc},
this paper presents two different ways of generalizing OGM and its development.
First, this paper specifies a parameterized family of 
algorithms that generalizes OGM,
and provides worst-case bounds 
on the function and gradient norm values
for this family.
Like the generalized forms of FGM~\cite{chambolle:15:otc,nesterov:05:smo}
being widely used and studied
(\eg,~\cite{attouch:18:fco,chambolle:15:otc,su:16:ade}),
we believe introducing the generalized OGM here
can be potentially useful.
Second, this paper optimizes 
the step coefficients of 
fixed-step first-order methods 
with respect to the rate of decrease of the cost function's gradient norm,
leading to a new algorithm called OGM-OG (OG for optimized over a gradient).
This development expands the choice of worst-case rate metrics for
optimizing first-order methods in~\cite{drori:14:pof,kim:16:ofo,kim:17:otc}
that focused on the cost function decrease leading to OGM.
We next briefly review the Performance Estimation Problem (PEP)~\cite{drori:14:pof} 
that was used in~\cite{drori:14:pof,kim:16:ofo,kim:17:otc} to develop OGM
and that we extensively use throughout the paper.

Drori and Teboulle~\cite{drori:14:pof}
cast a worst-case analysis
into an optimization problem
called PEP\footnote{
The original PEP was intractable to solve,
so a series of relaxation on the PEP was introduced in~\cite{drori:14:pof}
to make it possibly solvable,
which we review in Sec.~\ref{sec:pep,cost,relax}.}
\cite{drori:14:pof}
that examines the maximal absolute cost function inaccuracy
over all possible inputs (cost functions) to the optimization algorithm.
(See \eg,~\cite{deklerk:17:otw,drori:18:efo,drori:16:aov,kim:16:ofo,kim:17:otc,kim:18:ala,kim:18:ote-arxiv,lessard:16:aad,taylor:17:ewc-composite,taylor:17:ssc} for its extensions.)
Moreover, Drori and Teboulle~\cite{drori:14:pof}
optimized \emph{numerically} the step coefficients of first-order methods using PEP
for smooth convex minimization,
and found an algorithm
whose worst-case bound is lower than that of FGM,
but it required too much computation and memory to be appealing
for large-scale problems.
Building on their work,
the authors~\cite{kim:16:ofo,kim:17:otc}
found computationally and memory-wise efficient version,
called OGM,
and showed analytically that OGM satisfies
an analytical worst-case bound 
that is twice smaller than that of FGM.
Drori~\cite{drori:17:tei} showed that the OGM
is optimal 
for large-dimensional smooth convex minimization
over a general class of first-order methods
with either fixed or adaptive step sizes~\cite{drori:17:tei}.
This OGM has been numerically extended for nonsmooth composite convex problems
in~\cite{taylor:17:ewc-composite}.
In addition, this OGM-type algorithm was already studied
in the context of a proximal point method
\cite[Appendix]{guler:92:npp}.

Using the PEP approach~\cite{drori:14:pof},
this paper proposes
a generalized version of OGM (GOGM)
and analyzes its worst-case rates
in terms of \emph{both}
the decrease of the cost function
and the decrease of the norm
of the gradient of the cost function. 
The results complement the worst-case analysis 
of the OGM~\cite{kim:16:ofo,kim:17:otc},
and expands our understanding of OGM-type first-order methods.

This paper analyzes
the worst-case rate of the gradient norm
(in addition to that of the cost function)
because it
is important 
when dealing with dual problems,
considering that
the dual gradient norm
corresponds to the primal distance to feasibility
(see \eg,~\cite{devolder:12:dst,necoara:16:ica,nesterov:12:htm}).
While FGM has not been shown previously
to satisfy a rate $O(1/N^{1.5})$
for decreasing the gradient norm,
modified versions of FGM with such rate were studied
in~\cite{ghadimi:16:agm,monteiro:13:aah,nesterov:12:htm}.
This paper proves that FGM in fact does have that rate,
building upon~\cite{taylor:17:ssc}
that numerically conjectured such rate for FGM
using the gradient norm version of PEP.
For further acceleration of the worst-case gradient norm rate,
we optimize the step coefficients of first-order methods
with respect to the gradient norm using PEP
and propose an algorithm named OGM-OG
that belongs to the GOGM family
and has the best known analytical worst-case bound 
on rate of decrease of the gradient norm
among fixed-step first-order methods.

One can extend some aspects
of the approaches for generalizing OGM described in this paper
to other optimization algorithms and problems.
One direction we have already taken in~\cite{kim:18:ala}
aims to improve the fast iterative shrinkage/thresholding algorithm (FISTA)
\cite{beck:09:afi}
(that reduces to FGM for smooth convex problems)
for nonsmooth composite convex problems.
Naturally, this paper and~\cite{kim:18:ala}
use some similar approaches,
but they are different; 
the methods in~\cite{kim:18:ala}
when simplified to the smooth case correspond to a generalization of FGM
that differs from the GOGM. 
Another direction we have recently taken in
\cite{kim:18:ote-arxiv}
focuses on optimizing the step coefficients of first-order methods
with respect to the gradient norm
under the initial bounded function condition
that is different from the initial bounded distance condition
used in this paper.

Sec.~\ref{sec:prob}
defines the smooth convex problem
and the first-order methods.
Sec.~\ref{sec:conv} reviews and discusses worst-case analyses 
of a gradient method (GM), FGM, and OGM
for both the function value and the gradient norm.
Sec.~\ref{sec:conv}
also reviews first-order methods that guarantee an $O(1/N^{1.5})$ rate
for the gradient decrease.
Sec.~\ref{sec:pep,cost} reviews
the cost function form of PEP~\cite{drori:14:pof}
and reviews how the OGM~\cite{kim:16:ofo}
is derived using such PEP. 
Sec.~\ref{sec:pep,cost} then proposes a generalized version of OGM (GOGM)
using the cost function form of PEP,
and Sec.~\ref{sec:pep,sgnorm} 
provides a worst-case gradient norm bound for the GOGM
using the gradient form of PEP.
Then, Sec.~\ref{sec:pep,sgnorm} 
optimizes the step coefficients using the gradient form of PEP
and proposes the OGM-OG that belongs to the GOGM family.
Sec.~\ref{sec:pep,sgnorm}
also proves that FGM decreases the gradient norm
with a rate $O(1/N^{1.5})$. 
Sec.~\ref{sec:disc}
and Sec.~\ref{sec:conc} provide discussion and conclusion.

\section{Smooth convex problem and first-order methods} 
\label{sec:prob}

\subsection{Smooth convex problem}

We focus on the following smooth convex minimization problem
\begin{align}
\min_{\x\in\Reals^d} \;&\; f(\x)
\label{eq:prob} \tag{M}
,\end{align}
where the following additional conditions are assumed:
\begin{itemize}[leftmargin=40pt]
\item
$f\;:\;\Reals^d\rightarrow\Reals$
is a convex function of the type \cC, i.e., continuously differentiable with
Lipschitz continuous gradient:
\begin{align}
||\nabla f(\x) - \nabla f(\y)|| \le L||\x - \y||, \quad \forall \x, \y \in \Reals^d,
\end{align}
where $L > 0$ is the Lipschitz constant.
\item
The optimal set $X_*(f)=\argmin{\x\in\Reals^d} f(\x)$ is nonempty,
\ie, the problem~\eqref{eq:prob} is solvable.
\end{itemize}
We use \cF to denote the class of functions that satisfy the above conditions.
We also assume that
the distance between an initial point $\x_0$
and an optimal solution $\x_* \in X_∗(f)$ is bounded by
some $R > 0$, i.e., 
\begin{align}
||\x_0 - \x_*|| \le R
.\end{align}

\subsection{First-order methods}

To solve~\eqref{eq:prob}, 
we consider the following class of 
\emph{fixed-step} (or non-apdative-step)
first-order methods (\FO),
where the update step at $(i+1)$th iteration
is a weighted sum of the previous and current gradients
$\{\nabla f(\x_k)\}_{k=0}^i$
scaled by $\frac{1}{L}$
with fixed constant step coefficients $\{h_{i+1,k}\}_{k=0}^i$
that are not adaptive to the given $f$ and $\x_0$ (and thus $L$ and $R$).
This class~\FO includes GM, FGM, OGM,
and the methods proposed in this paper,
but excludes line-search-type methods.

\fbox{
\begin{minipage}[t]{0.85\linewidth}
\vspace{-10pt}
\begin{flalign}
&\quad \text{\bf Algorithm Class~\FO} & \nonumber \\
&\qquad \text{Input: } f\in \cF,\; \x_0\in\Reals^d. & \nonumber \\
&\qquad \text{For } i = 0,\ldots,N-1 & \nonumber \\
&\qquad \qquad \x_{i+1} = \x_i - \frac{1}{L}\sum_{k=0}^i \hkip \nabla f(\x_k). & \nonumber
\end{flalign}
\end{minipage}
} \vspace{5pt}

\section{Review of the worst-case analysis of~\FO}
\label{sec:conv}

This section reviews the worst-case analysis of existing~{\FO}s
(and simple variants thereof) 
in terms of bounds on the cost function and gradient norm.
Sec.~\ref{sec:conv_fv} reviews 
the worst-case cost function decrease of GM, FGM, and OGM.
Sec.~\ref{sec:grad} 
presents both~{\FO}s (including GM, FGM, OGM and some variants)
that have either $O(1/N)$ or $O(1/N^{1.5})$ rate
for the worst-case gradient decrease;
it also reviews an $O(1/N^2)$ lower bound of the worst-case rates of first-order methods
for decreasing the gradient norm.

\subsection{Function value worst-case analysis of~\FO}
\label{sec:conv_fv}

The simplest example of a~\FO is the following GM
that uses only the current gradient and the Lipschitz constant $L$
for the update.

\fbox{
\begin{minipage}[t]{0.85\linewidth}
\vspace{-10pt}
\begin{flalign}
&\quad \text{\bf Algorithm GM} & \nonumber \\
&\qquad \text{Input: } f\in \cF,\; \x_0\in\Reals^d. & \nonumber \\
&\qquad \text{For } i = 0,1,\ldots & \nonumber \\
&\qquad \qquad \x_{i+1} = \x_i - \frac{1}{L}\nabla f(\x_i) & \nonumber
\end{flalign}
\end{minipage}
} \vspace{5pt}

\noindent
This GM monotonically decreases the cost function~\cite{nesterov:04}
and satisfies the following tight\footnote{
A \emph{tight} worst-case bound denotes 
an inequality
where the equality holds for some function $f$.
For example,~\cite[Thm.~2]{drori:14:pof}
shows the bound~\eqref{eq:fv_gm} is tight.
}
worst-case bound~\cite[Thm.~1]{drori:14:pof}, for any $i\ge 0$,
\begin{align}
f(\x_i) - f(\x_*) \le \frac{LR^2}{4i + 2}
\label{eq:fv_gm}
.\end{align}

Among the class~\FO, the following two equivalent forms 
of FGM~\cite{nesterov:83:amf,nesterov:05:smo}
have been used widely
because they decrease the cost function
with the optimal rate $O(1/N^2)$.

\begin{center}
\fbox{
\small
\begin{minipage}[t]{0.45\textwidth}
\vspace{-10pt}
\begin{flalign}
&\text{\bf Algorithm FGM1} & \nonumber \\
&\text{Input: } f\in \cF,\; \x_0=\y_0\in\Reals^d, & \nonumber \\
&	\hspace{30pt} 
			t_0 = 1. & \nonumber \\
&\text{For } i = 0,1,\ldots & \nonumber \\
&\quad \y_{i+1} = \x_i - \frac{1}{L}\nabla f(\x_i) & \nonumber \\
&\quad t_{i+1} = \frac{1+\sqrt{1+4t_i^2}}{2}, & \nonumber \\
&\quad \x_{i+1} = \y_{i+1}
                + \frac{t_i - 1}{t_{i+1}}(\y_{i+1} - \y_i) \nonumber 
\end{flalign}
\end{minipage}\vline\hspace{8pt}%
\begin{minipage}[t]{0.45\textwidth}
\vspace{-10pt}
\begin{flalign}
&\text{\bf Algorithm FGM2} & \nonumber \\
&\text{Input: } f\in \cF,\; \x_0=\y_0\in\Reals^d, & \nonumber \\
&	\hspace{30pt} 
			t_0 = 1. & \nonumber \\
&\text{For } i = 0,1,\ldots & \nonumber \\
&\quad \y_{i+1} = \x_i - \frac{1}{L}\nabla f(\x_i) & \nonumber \\
&\quad \z_{i+1} = \x_0 - \frac{1}{L}\sum_{k=0}^it_k \nabla f(\x_k) & \nonumber \\
&\quad t_{i+1} = \frac{1+\sqrt{1+4t_i^2}}{2}, & \nonumber \\
&\quad \x_{i+1} = \paren{1 - \frac{1}{t_{i+1}}}\y_{i+1}
                + \frac{1}{t_{i+1}}\z_{i+1} & \nonumber
\end{flalign}
\end{minipage}
} \vspace{5pt}
\end{center}

\noindent
Specifically, the FGM1 and FGM2 iterates satisfy 
the following worst-case cost function bounds
\cite{kim:16:ofo,nesterov:83:amf,nesterov:05:smo} for any $i\ge 1$:
\begin{align}
f(\y_i) - f(\x_*) \le \frac{LR^2}{2t_{i-1}^2}
	\le \frac{2LR^2}{(i+1)^2},
\quad\text{and}\quad
f(\x_i) - f(\x_*) \le \frac{LR^2}{2t_i^2}
	\le \frac{2LR^2}{(i+2)^2}
\label{eq:fv_fgm}
,\end{align}
where the parameter $t_i$ satisfies 
\begin{align}
t_i^2 = \sum_{l=0}^i t_l
\quad\text{and}\quad 
t_i \ge \frac{i+2}{2}
\quad\text{ for all } i.
\label{eq:t_rule}
\end{align}
A generalized form of FGM in~\cite{chambolle:15:otc} 
uses parameters $t_i$
satisfying $t_0=1$ and $t_i^2 \le t_{i-1}^2 + t_i$,
including the choice $t_i = \frac{i+a}{a}$ for any $a\ge2$.
There is another generalized form of FGM in~\cite{nesterov:05:smo},
and these generalized forms of FGM
have been widely used and studied
(\eg,~\cite{attouch:18:fco,chambolle:15:otc,su:16:ade}).
Similarly this paper studies generalizations of
the OGM.

Building upon~\cite{drori:14:pof}
that optimized numerically the step coefficients over the cost function form of PEP, 
the authors~\cite{kim:16:ofo} developed
the following two equivalent forms of OGM,
as reviewed in Sec.~\ref{sec:pep,cost}.

\begin{center}
\fbox{
\small
\begin{minipage}[t]{0.45\textwidth}
\vspace{-10pt}
\begin{flalign*}
&\text{\bf Algorithm OGM1} & \\
&\text{Input: } f\in \cF,\; \x_0=\y_0\in\Reals^d, & \\
&	\hspace{30pt} \theta_0 = 1. & \\
&\text{For } i = 0,\ldots,N-1 & \\
&\quad \y_{i+1} = \x_i - \frac{1}{L}\nabla f(\x_i) & \\
&\quad \theta_{i+1}
                = \begin{cases}
                        \frac{1+\sqrt{1+4\theta_i^2}}{2}, & i \le N-2 \\
                        \frac{1+\sqrt{1+8\theta_i^2}}{2}, & i = N-1
                        \end{cases} & \\
&\quad \x_{i+1} = \y_{i+1}
                + \frac{\theta_i - 1}{\theta_{i+1}}(\y_{i+1} - \y_i) & \\
&\hspace{70pt}
                + \frac{\theta_i}{\theta_{i+1}}(\y_{i+1} - \x_i) &
\end{flalign*}
\end{minipage}\vline\hspace{8pt}%
\begin{minipage}[t]{0.45\textwidth}
\vspace{-10pt}
\begin{flalign*}
&\text{\bf Algorithm OGM2} & \\
&\text{Input: } f\in \cF,\; \x_0=\y_0\in\Reals^d, & \\
&	\hspace{30pt} \theta_0 = 1. & \\
&\text{For } i = 0,\ldots,N-1 & \\
&\quad \y_{i+1} = \x_i - \frac{1}{L}\nabla f(\x_i) & \\
&\quad \z_{i+1} = \x_0 - \frac{1}{L}\sum_{k=0}^i2\theta_k \nabla f(\x_k) & \\
&\quad \theta_{i+1}
                = \begin{cases}
                        \frac{1+\sqrt{1+4\theta_i^2}}{2}, & i \le N-2 \\
                        \frac{1+\sqrt{1+8\theta_i^2}}{2}, & i = N-1
                        \end{cases} & \\
&\quad \x_{i+1} = \paren{1 - \frac{1}{\theta_{i+1}}}\y_{i+1}
                + \frac{1}{\theta_{i+1}}\z_{i+1} &
\end{flalign*}
\end{minipage}
} \vspace{5pt}
\end{center}

\noindent
The OGM iterates
satisfy the following worst-case 
cost function bounds~\cite{kim:16:ofo,kim:17:otc}:
\begin{align}
f(\y_i) - f(\x_*) &\le \frac{LR^2}{4\theta_{i-1}^2}
	\le \frac{LR^2}{(i+1)^2},
\label{eq:fv_ogm_}
\end{align}
for any $1\le i\le N$, and
\begin{align}
f(\x_N) - f(\x_*) &\le \frac{LR^2}{2\theta_N^2} 
	\le \frac{LR^2}{(N+1)(N+1+\sqrt{2})}
\label{eq:fv_ogm}
.\end{align}
The parameter sequence $\theta_i$ satisfies 
\begin{align}
\theta_i^2 = \begin{cases}
\sum_{l=0}^{i+1}\theta_l, & i \le N-1, \\
2\sum_{l=0}^{N-1}\theta_l + \theta_N, & i = N, 
\end{cases}
\quad\text{and}\quad 
\theta_i \ge \begin{cases}
\frac{i+2}{2}, & i \le N-1, \\
\frac{i+1}{\sqrt{2}}, & i = N,
\end{cases}
\label{eq:theta_rule}
\end{align}
which is equivalent to $t_i$~\eqref{eq:t_rule} 
except at the final iteration.
The bounds~\eqref{eq:fv_ogm_} and~\eqref{eq:fv_ogm} of OGM
are about twice smaller than the bounds~\eqref{eq:fv_fgm} of FGM,
so OGM decreases the cost function faster than FGM
in the worst case
(and often in practice~\cite{kim:15:aof}).
In addition, the bound~\eqref{eq:fv_ogm} on the final iterate $\x_N$
is tight and
satisfies the optimal worst-case bound of 
general first-order methods
including both~\FO and adaptive-step first-order methods,
when the condition $d\ge N+1$ holds~\cite{drori:17:tei}.

The additional term $\frac{\theta_i}{\theta_{i+1}}(\y_{i+1} - \x_i)$ of OGM1
and the additional constant $2$ for the update of $\z_i$ of OGM2,
compared to FGM1 and FGM2 respectively,
along with the parameter $\theta_N$,
are what make
OGM optimal (for $d \ge N+1$).
One of the main goals of this paper is to
generalize the form of OGM
and analyze the worst-case rate of such generalized OGM
in terms of both the function value and the gradient norm,
complementing the bounds~\eqref{eq:fv_ogm_} and~\eqref{eq:fv_ogm} 
on the function value of OGM. 

The next section studies the worst-case rate of the gradient of~\FO.

\subsection{Gradient norm worst-case analysis of~\FO}
\label{sec:grad}

When tackling dual problems,
it is known that the gradient norm worst-case rate is important 
in addition to the function value worst-case rate
because the dual gradient norm
is related to the primal distance to feasibility
(see \eg, \cite{devolder:12:dst,necoara:16:ica,nesterov:12:htm}).
One simple way to find a (loose) worst-case bound for the gradient norm
is to use the well-known convex inequality
for convex functions with $L$-Lipschitz continuous gradients~\cite{nesterov:04}:
\begin{align}
\frac{1}{2L}||\nabla f(\x)||^2 
	\le f(\x) - f\paren{\x - \frac{1}{L}\nabla f(\x)}
	\le f(\x) - f(\x_*), \quad \forall \x\in\Reals^d
\label{eq:g_upper}
,\end{align}
as discussed in~\cite{nesterov:04,taylor:17:ssc}.
Combining the bounds~\eqref{eq:fv_gm},~\eqref{eq:fv_fgm}
and the inequality~\eqref{eq:g_upper},
for any $i\ge1$,
the GM iterates satisfy
\begin{align}
||\nabla f(\x_i)|| \le \sqrt{2L(f(\x_i) - f(\x_*))} \le \frac{LR}{\sqrt{2i+1}}
\label{eq:gg_gm}
,\end{align}
and the iterates of FGM satisfy
\begin{align}
||\nabla f(\y_i)|| \le \frac{LR}{t_{i-1}} \le \frac{2LR}{i+1},
\quad\text{and}\quad
||\nabla f(\x_i)|| \le \frac{LR}{t_i} \le \frac{2LR}{i+2}
\label{eq:g_fgm}
.\end{align}
Similarly for any $1\le i\le N$, 
the OGM iterates with the bounds~\eqref{eq:fv_ogm_} and~\eqref{eq:fv_ogm}
satisfy
\begin{align}
||\nabla f(\y_i)|| \le \frac{LR}{\sqrt{2}\theta_{i-1}} \le \frac{\sqrt{2}LR}{i+1},
\quad\text{and}\quad
||\nabla f(\x_N)|| \le \frac{LR}{\theta_N} \le \frac{\sqrt{2}LR}{N+1}
\label{eq:g_ogm}
.\end{align}
Unfortunately, using the inequality~\eqref{eq:g_upper} 
provides at best an $O(1/N)$ bound
due to the optimal rate $O(1/N^2)$ of the function decrease.
Furthermore, in general using~\eqref{eq:g_upper}
need not lead to tight worst-case bounds
on the gradient norm.

Using a different approach, 
a smaller $O(1/N)$ worst-case bound for the gradient norm of GM
was derived in~\cite{nesterov:12:htm},
as reviewed in next section.
While the bounds~\eqref{eq:gg_gm},~\eqref{eq:g_fgm}, and~\eqref{eq:g_ogm}
are not guaranteed to be tight,
the next section shows that
the worst-case gradient bound~\eqref{eq:g_ogm} 
on the final iterate $\x_N$ of OGM 
is in fact tight
and thus has the same disappointingly slow $O(1/N)$ worst-case bound
on the gradient norm as GM.

\subsubsection{\FO with rate $O(1/N)$ for decreasing the gradient norm}

This section uses the following lemma
stating that GM monotonically decreases the gradient.

\begin{lemma}
\label{lem:gmono}
\cite[Lemma~2.4]{necoara:16:ica}
The GM monotonically decreases the gradient norm, i.e.,
\begin{align}
\bigg|\bigg|\nabla f\paren{\x - \frac{1}{L}\nabla f(\x)}\bigg|\bigg| \le ||\nabla f(\x)||
\label{eq:g_desc}
.\end{align}
\end{lemma}

The following theorem reviews
a simple proof in~\cite{nesterov:12:htm} 
that provides a worst-case gradient norm bound for GM 
with rate $O(1/N)$
that is smaller than~\eqref{eq:gg_gm},
where~\cite[Thm.~6.1]{necoara:16:ica}
additionally considers Lemma~\ref{lem:gmono}.
\begin{theorem}
\cite[Thm.~6.1]{necoara:16:ica},~\cite{nesterov:12:htm}
Let $f\;:\;\Reals^d\rightarrow\Reals$ be $\cF$ 
and let $\x_0,\cdots,\x_N\allowbreak\in\Reals^d$ be
generated by GM. Then for any $N\ge1$,
\begin{align}
\min_{i\in\{0,\ldots,N\}} ||\nabla f(\x_i)|| = ||\nabla f(\x_N)|| 
	\le \frac{\sqrt{2}LR}{\sqrt{N(N+2)}}
\label{eq:g_gm}
.\end{align}
\begin{proof}
Lemma~\ref{lem:gmono}
implies the first equality in~\eqref{eq:g_gm}.
Using~\eqref{eq:fv_gm},~\eqref{eq:g_upper},~\eqref{eq:g_desc} yields:
\begin{align*}
\frac{LR^2}{4m+2} 
	&\stackrel{\eqref{eq:fv_gm}}{\ge} f(\x_m) - f(\x_*) 
        \stackrel{\eqref{eq:g_upper}}{\ge} 
		f(\x_{N+1}) - f(\x_*) + \frac{1}{2L}\sum_{i=m}^N||\nabla f(\x_i)||^2 \\
	&\stackrel{\eqref{eq:g_desc}}{\ge} \frac{N-m+1}{2L} ||\nabla f(\x_N)||^2
,\end{align*}
which is equivalent to~\eqref{eq:g_gm}
using $m = \fNh$ for which $m \ge \frac{N-1}{2}$ and $N-m \ge \frac{N}{2}$.
\end{proof}
\end{theorem}

Inspired by the conjecture in~\cite[Sec.~4.1.3]{taylor:17:ssc},
the following theorem shows that
the $O(1/N)$ rate of the worst-case gradient norm bound~\eqref{eq:g_gm} of GM
is tight up to a constant.

\begin{theorem}
Let $\x_0,\cdots,\x_N \in \Reals^d$ be generated by
GM. Then for any $N\ge1$,
\begin{align}
\frac{LR}{N+1} 
\le \max_{\substack{f\in\cF, \\ \x_* \in X_∗(f), \\ ||\x_0 - \x_*||\le R}}\min_{i\in\{0,\ldots,N\}} ||\nabla f(\x_i)||
= \max_{\substack{f\in\cF, \\ \x_* \in X_∗(f), \\ ||\x_0 - \x_*||\le R}}||\nabla f(\x_N)||
\label{eq:g_gm_low}
,\end{align}
where the 
inequality in~\eqref{eq:g_gm_low} is achieved
by the following function in \cF:
\begin{align}
\psi(\x) = \begin{cases}
        \frac{LR}{N+1}||\x|| - \frac{LR^2}{2(N+1)^2}, & ||\x|| \ge \frac{R}{N+1}, \\
        \frac{L}{2}||\x||^2, & ||\x|| < \frac{R}{N+1}.
        \end{cases}
\end{align}
\begin{proof}
Lemma~\ref{lem:gmono}
implies the equality in~\eqref{eq:g_gm_low}.
Starting from $\x_0 = R\bmnu$, where $\bmnu$ is a unit vector,
the GM iterates are
\begin{align*}
\x_i = \paren{1 - \frac{i}{N+1}}R\bmnu, \quad
\nabla \psi(\x_i) = \frac{LR}{N+1}\bmnu, \quad i=0,\ldots,N,
\end{align*}
which implies the inequality~\eqref{eq:g_gm_low}.
\end{proof}
\end{theorem}

We next show that the bound~\eqref{eq:g_ogm} 
for the gradient norm at the final iterate $\x_N$ of OGM is tight
and and that its worst-case function is a simple quadratic function.
Note that OGM was derived by optimizing a worst-case bound 
on the cost function decrease
and its behavior in terms of gradient norms
was not investigated previously.

\begin{theorem}
\label{thm:g_ogm_tight}
Let $\x_0,\cdots,\x_N\in\Reals^d$ be
generated by OGM. Then for any $N\ge1$,
\begin{align}
	\max_{\substack{f\in\cF, \\ \x_* \in X_∗(f), \\ ||\x_0 - \x_*||\le R}} \min_{i\in\{0,\ldots,N\}} 
		||\nabla f(\x_i)|| 
		= \max_{\substack{f\in\cF, \\ \x_* \in X_∗(f), \\ ||\x_0 - \x_*||\le R}} ||\nabla f(\x_N)||
	= \frac{LR}{\theta_N}
	\;\bigg(\! \le \frac{\sqrt{2}LR}{N+1}\bigg)
\label{eq:g_ogm_tight}
,\end{align}
where the worst-case function in \cF for OGM in terms of the gradient norm is
the quadratic function 
$\phi(\x) = \frac{L}{2}||\x||^2$.
\begin{proof}
See Appendix~\ref{appen1}.
\end{proof}
\end{theorem}

Comparing~\eqref{eq:g_gm} and~\eqref{eq:g_ogm_tight},
we see that GM and OGM 
have essentially similar
worst-case gradient norm bounds.
This is a dilemma
because OGM is the fastest~\FO in terms of the worst-case cost function bound,
but is as slow as GM in terms of the worst-case gradient norm bound.
Therefore, one of the main goals of this paper 
is to study optimizing the step coefficients of~\FO
using PEP with respect to the gradient norm
in Sec.~\ref{sec:opt,grad,pep}.

We next discuss
the specific~\FO in~\cite{nesterov:12:htm}
that decreases the gradient norm
with a faster $O(1/N^{1.5})$ rate.

\subsubsection{\FO with rate $O(1/N^{1.5})$ 
for decreasing the gradient norm}

Searching for a~\FO that decreases the gradient norm faster 
than the $O(1/N)$ rate of GM (and OGM),
Nesterov~\cite{nesterov:12:htm} 
(among other variants of FGM~\cite{ghadimi:16:agm,monteiro:13:aah})
considered performing
FGM for the first \m iterations,
and GM for the remaining iterations.
He showed that this method,
which we denote FGM-\m, 
satisfies a fast rate $O(1/N^{1.5})$
for decreasing the gradient norm.
In~\cite{devolder:12:dst,necoara:16:ica,nesterov:12:htm},
FGM-\m for $m=\fNh$ was used to solve dual problems.
To pursue a faster worst-case rate (in terms of the constant factor),
we consider here another variant
that performs OGM for the first \m iterations
and GM for the remaining iterations,
which we denote OGM-\m. 

\noindent
\hspace{5pt}
\fbox{
\small
\begin{minipage}[t]{0.85\linewidth}
\vspace{-10pt}
\begin{flalign*}
&\quad \text{\bf Algorithm OGM-\m} & \\
&\qquad \text{Input: } f\in \cF,\; \x_0=\y_0\in\Reals^d,\;\vartheta_0 = 1,\;
	m\in\{1,\ldots,N-1\}. & \\
&\qquad \text{For } i = 0,\ldots,\m - 1 & \\
&\qquad \qquad \y_{i+1} = \x_i - \frac{1}{L}\nabla f(\x_i) & \\
&\qquad \qquad \vartheta_{i+1} = \begin{cases}
                        \frac{1+\sqrt{1+4\vartheta_i^2}}{2}, & i\le\m - 2 \\
                        \frac{1+\sqrt{1+8\vartheta_i^2}}{2}, & i=\m - 1 \\
                        \end{cases}
                & \\
&\qquad \qquad \x_{i+1}
                = 
                \y_{i+1} + \frac{\vartheta_i-1}{\vartheta_{i+1}}(\y_{i+1} - \y_i) 
			+ \frac{\vartheta_i}{\vartheta_{i+1}}(\y_{i+1} - \x_i)
		& \\
&\qquad \text{For } i = \m,\ldots,N-1 & \\
&\qquad \qquad \x_{i+1} = \x_i - \frac{1}{L}\nabla f(\x_i) &
\end{flalign*}
\end{minipage}
} 

\noindent
The following theorem bounds the gradient norm of the OGM-\m iterates,
inspired by the proof in~\cite{necoara:16:ica,nesterov:12:htm}
for the worst-case gradient norm bound of the FGM-\m iterates.
The worst-case bound of FGM-\m 
in~\cite{necoara:16:ica,nesterov:12:htm}
is asymptotically $\sqrt{2}$-times larger than
the following new bound~\eqref{eq:g_ogmh} 
for OGM-\m. 

\begin{theorem}
\label{thm:ogmh}
Let $f\;:\;\Reals^d\rightarrow\Reals$ be $\cF$
and let $\x_0,\cdots\x_N\in\Reals^d$ be
generated by OGM-\m for $1 \le m \le N-1$. Then for any $N\ge1$,
\begin{align}
\min_{i\in\{0,\ldots,N\}} ||\nabla f(\x_i)||
        \le ||\nabla f(\x_N)||
	\le \frac{\sqrt{2}LR}{(m+1)\sqrt{N-m+1}}
\label{eq:g_ogmh}
.\end{align}
\begin{proof}
Using~\eqref{eq:fv_ogm},~\eqref{eq:g_upper},~\eqref{eq:g_desc} yields:
\begin{align*}
\frac{LR^2}{2\vartheta_m^2}
        &\stackrel{\eqref{eq:fv_ogm}}{\ge} f(\x_m) - f(\x_*)
        \stackrel{\eqref{eq:g_upper}}{\ge}
                f(\x_{N+1}) - f(\x_*) + \frac{1}{2L}\sum_{i=m}^N||\nabla f(\x_i)||^2 \\
	&\stackrel{\eqref{eq:g_desc}}{\ge} \frac{N-m+1}{2L}||\nabla f(\x_N)||^2
,\end{align*}
which is equivalent to~\eqref{eq:g_ogmh}
using 
$\vartheta_m \ge \frac{m+1}{\sqrt{2}}$
that is implied by~\eqref{eq:theta_rule}.
\end{proof}
\end{theorem}
The bound~\eqref{eq:g_ogmh}
is minimized
at a point close to $m=\fNt$,
leading to its (approximately) smallest constant $\frac{3\sqrt{6}}{2}$
with the rate $O(1/N^{1.5})$. 

Other variants of FGM 
having $O(1/N^{1.5})$ worst-case gradient bounds
were derived in~\cite{ghadimi:16:agm,monteiro:13:aah}.
Such variations of FGM (including FGM-$m$) were derived
since, prior to this paper, it was unknown 
whether or not FGM 
decreases the gradient norm with the rate $O(1/N^{1.5})$;
this rate for the gradient norm of FGM
was conjectured numerically in~\cite{taylor:17:ssc}.
Sec.~\ref{sec:pep,sgnorm,fgm} below
uses the PEP to show for the first time
the rate $O(1/N^{1.5})$ for the gradient decrease of the FGM.
The bound~\eqref{eq:g_ogmh} of OGM-\m for decreasing the gradient
is smaller than the bounds for the FGM variants
in~\cite{ghadimi:16:agm,monteiro:13:aah},
and
Sec.~\ref{sec:opt,grad,pep} below shows that our proposed methods
have worst-case bounds even lower than~\eqref{eq:g_ogmh}.

The preceding sections have focused on tight or upper worst-case bounds
of the gradient norm decrease of first-order methods,
whereas the next section
reviews a lower bound for the worst-case gradient norm decrease in
\cite{nemirovsky:92:ibc},
illustrating the best achievable worst-case rate
of the gradient norm decrease
for any first-order method 
(with either fixed-step or adaptive-step approaches).

\subsubsection{A lower bound of the worst-case rates of 
first-order methods
for decreasing the gradient norm}
\label{sec:low,grad}

For completeness, this section reviews
a lower bound on the worst-case rate of any first-order method
in terms of the gradient norm values
for smooth convex \emph{quadratic} functions~\cite{nemirovsky:92:ibc}.
Lower bounds on the function value
were studied
for convex \emph{quadratic} functions in~\cite{nemirovsky:92:ibc},
and for smooth convex functions
in~\cite{drori:17:tei,nesterov:04}.

When the condition $d\ge 2N+3$ holds,
a worst-case gradient norm bound 
of any first-order method generating $\x_N$ after $N$ iterations
has rate $O(1/N^2)$ at best, for convex quadratic $f$,
\ie,
has the following lower bound~\cite[Sec.~2.3.B]{nemirovsky:92:ibc}:
\begin{align}
\frac{LR}{4e^2(N+1)^2} 
	\le \max_{\substack{f\in\cQ, \\ \x_* \in X_∗(f), \\ ||\x_0 - \x_*||\le R}}||\nabla f(\x_N)|| 
\label{eq:low}
,\end{align}
where 
$
\cQ := \left\{f\;:\;
        f(\x) \equiv \frac{1}{2}\x^\top\Q\x + \p^\top\x + \rr
                \text{ for } \x\in\reals^d,\; 
        \Q\succeq\Zero,\; ||\Q||\le L
	\right\}
.$
Since $\cQ\subset\cF$, the lower bound~\eqref{eq:low}
for convex quadratic functions 
also applies to smooth convex functions.

A regularization technique in~\cite{nesterov:12:htm}
achieves the rate $O(1/N^2)$ up to a logarithmic factor.
However, its adaptive step coefficients
require knowing $R$ in advance which is undesirable in practice.
To our knowledge, 
whether there exists any \FO satisfying such rate
is an open question.
Instead, this paper discusses a way to develop \FO
that achieves an $O(1/N^{1.5})$ gradient norm bound 
with the smallest constant
among known \FO.

\section{Relaxation and optimization of the cost function form of PEP}
\label{sec:pep,cost}

This section reviews
a relaxation of the cost function form of PEP~\cite{drori:14:pof} and
reviews how~\cite{drori:14:pof,kim:16:ofo,kim:17:otc}
optimized the step coefficients of the~\FO class
over the cost function form of PEP,
leading to OGM.
Then, we 
propose a parameterized family of algorithms
that generalizes OGM,
and analyze the worst-case cost function decrease of the generalized OGM family.

\subsection{Review: Relaxation for the cost function form of PEP}
\label{sec:pep,cost,relax}

The worst-case bound on the cost function
for a~\FO having given step coefficients $\bmh:=\{h_{i+1,k}\}$
corresponds to a solution of the following PEP problem~\cite[Prob.~(P)]{drori:14:pof}:
\begin{align}\mathcal{B}_{\mathrm{P}}(\bmh,N,d,L,R) :=\;
& \max_{\substack{f\in\cF, \\
	\x_0,\cdots,\x_N\in\Reals^d, \\
	\x_*\in X_*(f), \\
	||\x_0 - \x_*|| \le R}}        
	f(\x_N) - f(\x_*)
        \label{eq:PEP} \tag{P} \\
        &\st \; \x_{i+1} = \x_i - \frac{1}{L}\sum_{k=0}^i \hkip \nabla f(\x_k),
                \quad i=0,\ldots,N-1. \nonumber
\end{align}
Since problem~\eqref{eq:PEP} is impractical to solve
due to its functional constraint 
$f\in\cF$,
\cite{drori:14:pof} relaxed it by 
the following finite set of inequalities
satisfied by $f$~\cite[Thm.~2.1.5]{nesterov:04}:
\begin{align}
\frac{1}{2L}||\nabla f(\x_i) - \nabla f(\x_j)||^2
	\le f(\x_i) - f(\x_j) - \inprod{\nabla f(\x_j)}{\x_i - \x_j}
\label{eq:ineq}
\end{align}
for $i,j=0,1,\ldots,N,*$.
Then, a matrix
$\G = [\g_0,\cdots,\g_N]^\top \in \Reals^{(N+1)\times d}$
and a vector $\del = [\delta_0,\cdots,\delta_N]^\top \in \Reals^{N+1}$
with
\begin{align*}
\g_i := \frac{1}{L||\x_0 - \x_*||}\nabla f(\x_i),
\text{ and }
\delta_i := \frac{1}{L||\x_0 - \x_*||^2}(f(\x_i) - f(\x_*))
\end{align*}
are introduced
to represent gradient vectors and function values respectively 
in the set of~\eqref{eq:ineq}.
This leads to
a finite-dimensional relaxation of problem~\eqref{eq:PEP}
\cite[Prob.~(Q)]{drori:14:pof}:
\begin{align}\mathcal{B}_{\mathrm{P1}}(\bmh,N,d,L,R) :=\;
& \max_{\substack{\G\in\Reals^{(N+1)\times d}, \\
        \del\in\Reals^{N+1}}}        
	LR^2\delta_N
        \label{eq:relPEP} \tag{P1} \\
        &\st \; \Tr{\G^\top\A_{i,j}(\bmh)\G} \le \delta_i - \delta_j,
                \quad i<j=0,\ldots,N, \nonumber  \\
        &\quad\;\;\; \Tr{\G^\top\B_{i,j}(\bmh)\G} \le \delta_i - \delta_j,
		\quad j<i=0,\ldots,N, \nonumber \\
	&\quad\;\;\; \Tr{\G^\top\C_i\G} \le \delta_i, 
		\quad i=0,\ldots,N, \nonumber \\
	&\quad\;\;\; \Tr{\G^\top\D_i(\bmh)\G + \bmnu\u_i^\top\G} \le - \delta_i,
		\quad i=0,\ldots,N, \nonumber
\end{align}
for any given unit vector $\bmnu\in\Reals^d$,
where $\u_i = \e_{i+1} \in\Reals^{N+1}$
is the $(i+1)$th standard basis vector.
Note that
$\Tr{\G^\top\u_i\u_j^\top\G} = \inprod{\g_i}{\g_j}$ by definition.
The matrices $\A_{i,j}(\bmh), \B_{i,j}(\bmh), \C_i, \D_i(\bmh)$ are defined as
\begin{align}
\begin{cases}
\A_{i,j}(\bmh) :=
        \frac{1}{2}(\u_i - \u_j)(\u_i - \u_j)^\top
        + \frac{1}{2}\sum_{l=i+1}^j\sum_{k=0}^{l-1}
        \hkl (\u_j\u_k^\top + \u_k\u_j^\top), & \\
\B_{i,j}(\bmh) :=
	\frac{1}{2}(\u_i - \u_j)(\u_i - \u_j)^\top
        - \frac{1}{2}\sum_{l=j+1}^i\sum_{k=0}^{l-1}
        \hkl (\u_j\u_k^\top + \u_k\u_j^\top), & \\
\C_i := \frac{1}{2}\u_i\u_i^\top, & \\	
\D_i(\bmh) := 
	\frac{1}{2}\u_i\u_i^\top + \frac{1}{2}\sum_{j=1}^i\sum_{k=0}^{j-1}
        \hkj(\u_i\u_k^\top + \u_k\u_i^\top). &
\end{cases}
\label{eq:ABCDF}
\end{align}

In~\cite[Prob.~(Q$'$)]{drori:14:pof},
problem~\eqref{eq:relPEP} is further relaxed
by discarding some constraints to yield
\begin{align}\mathcal{B}_{\mathrm{P2}}(\bmh,N,d,L,R) :=\;
& \max_{\substack{\G\in\Reals^{(N+1)\times d}, \\
        \del\in\Reals^{N+1}}}
        LR^2\delta_N
        \label{eq:pPEP} \tag{P2} \\
        &\st \; \Tr{\G^\top\A_{i-1,i}(\bmh)\G} \le \delta_{i-1} - \delta_i,
                \quad i=1,\ldots,N, \nonumber  \\
        &\quad\;\;\; \Tr{\G^\top\D_i(\bmh)\G + \bmnu\u_i^\top\G} \le - \delta_i,      
                \quad i=0,\ldots,N, \nonumber
\end{align}
for any given unit vector $\bmnu\in\Reals^d$.
We explicitly illustrate the relaxation 
from~\eqref{eq:relPEP} to~\eqref{eq:pPEP}
because Sec.~\ref{sec:pep,sgnorm} 
uses a similar but different relaxation.
Taylor~\etal~\cite{taylor:17:ssc} avoided this step
to analyze a tight worst-case bound of~\eqref{eq:PEP}
(under a large-scale condition $d\ge N+2$
\cite[Thm.~5]{taylor:17:ssc});
however, this relaxation~\eqref{eq:pPEP}
facilitates the analysis in~\cite{drori:14:pof,kim:16:ofo,kim:17:otc} and in this paper.

Replacing $\max_{\G,\del} LR^2\delta_N$
by $\min_{\G,\del} \{-\delta_N\}$ for convenience in~\eqref{eq:pPEP},
the Lagrangian of the corresponding constrained minimization problem
with dual variables $\bmlam = (\lambda_1,\cdots,\lambda_N)^\top \in\Reals^N$ 
and $\bmtau = (\tau_0,\cdots,\tau_N)^\top\in\Reals^{N+1}$ 
for the first and second constraint inequalities of~\eqref{eq:pPEP} 
respectively becomes
\begin{align}
\cL(\G,\del,\bmlam,\bmtau;\bmh)
        =& - \delta_N + \sum_{i=1}^N\lambda_i(\delta_i - \delta_{i-1})
                + \sum_{i=0}^N\tau_i\delta_i \\
	&\qquad\qquad\qquad\qquad + \Tr{\G^\top\S(\bmh,\bmlam,\bmtau)\G + \bmnu\bmtau^\top\G},
	\nonumber
\end{align}
where
\begin{align}
\S(\bmh,\bmlam,\bmtau) := \sum_{i=1}^N \lambda_i\A_{i-1,i}(\bmh)
                + \sum_{i=0}^N\tau_i\D_i(\bmh).
\label{eq:S}
\end{align}
Then, we have the following dual problem of~\eqref{eq:pPEP} 
that one could use to compute a valid upper bound of~\eqref{eq:PEP} 
using a semidefinite program (SDP)
for given $\bmh$
\cite[Prob.~(DQ$'$)]{drori:14:pof}:
\begin{align}
\mathcal{B}_{\mathrm{D}}(\bmh,N,L,R) :=\;
&\min_{\substack{(\bmlam,\bmtau)\in\Lambda, \\ \gamma\in\Reals}}
        \left\{
        \frac{1}{2}LR^2\gamma\;:\;
        \left(\begin{array}{cc}
                \S(\bmh,\bmlam,\bmtau) & \frac{1}{2}\bmtau \\
                \frac{1}{2}\bmtau^\top & \frac{1}{2}\gamma
        \end{array}\right)
        \succeq\Zero
        \right\}
        \label{eq:D} \tag{D}
,\end{align}
where
\begin{align}
\Lambda = \left\{(\bmlam,\bmtau)\in\Reals_+^{2N+1}
                \;:\;
                \begin{array}{l}
                \tau_0 = \lambda_1,\;\; \lambda_N + \tau_N = 1, \\
                \lambda_i - \lambda_{i+1} + \tau_i = 0, \; i=1,\ldots,N-1,
                \end{array}
                \right\}
\label{eq:Lam}
.\end{align}
The next section reviews the analytical solution to this upper bound~\eqref{eq:D} 
for OGM, 
instead of using a numerical SDP solver.

\subsection{Review: Optimizing step coefficients 
for the cost function form of PEP}
\label{sec:pep,cost,opt}

Drori and Teboulle~\cite{drori:14:pof} optimized numerically
the step coefficients $\bmh$ over 
the simple SDP problem~\eqref{eq:D} as follows
\cite[Prob.~(BIL)]{drori:14:pof}:
\begin{align}
\hat{\bmh}_{\mathrm{D}}
        := 
	\argmin{\bmh\in\Reals^{N(N+1)/2}} \mathcal{B}_{\mathrm{D}}(\bmh,N,L,R)
	\tag{HD}
\label{eq:HD}
.\end{align}
The problem~\eqref{eq:HD} is bilinear,
and \cite[Thm.~3]{drori:14:pof}\footnote{
\cite[Thm.~3]{drori:14:pof} has typos
that are fixed in~\cite[Eq. (6.3)]{kim:16:ofo}.}
used a convex relaxation technique
to make it solvable by numerical methods.

In~\cite[Lemma~4]{kim:16:ofo},
we solved~\eqref{eq:HD} analytically
yielding the optimized step coefficients
\begin{align}
\hkip &= \begin{cases}
                \frac{1}{\theta_{i+1}}
                        \paren{2\theta_k - \sum_{j=k+1}^i\hkj}, & k=0,\ldots,i-1, \\
                1 + \frac{2\theta_i - 1}{\theta_{i+1}}, & k=i,
        \end{cases}
       \label{eq:hh_ogm}
\end{align}
for $\theta_i$ in~\eqref{eq:theta_rule}.
Fortuitously, the optimized coefficients~\eqref{eq:hh_ogm}
lead to equivalent computationally 
efficient OGM1 and OGM2 forms~\cite[Prop.~3, 4 and 5]{kim:16:ofo},
and the bound~\eqref{eq:fv_ogm}  
for the final secondary iterate $\x_N$ of
OGM is implied by~\cite[Lemma~4]{kim:16:ofo}.
Recently,
Drori~\cite{drori:17:tei}
showed that the OGM is optimal for $d \ge N+1$,
implying that
optimizing over the relaxed bound~\eqref{eq:D} in~\eqref{eq:HD}
for simplicity
is equivalent to optimizing 
over the exact worst-case cost function bound~\eqref{eq:PEP}
when $d\ge N+1$.

One could use a SDP solver to compute a numerical bound
from~\eqref{eq:D} 
for any~\FO; however,
deriving an analytical bound using~\eqref{eq:D} is difficult
for the primary sequence $\{\y_i\}$ of OGM.
Therefore, we devised a new relaxed bound in~\cite{kim:17:otc}
similar to~\eqref{eq:D},
which we review next. 

\subsection{Review: Another cost function form of relaxed PEP for the primary sequence of OGM}

An upper bound of the worst-case bound on $f(\y_{N+1}) - f(\x_*)$
for~\FO with step coefficients \bmh
and $\y_{N+1} = \x_N - \frac{1}{L}f(\x_N)$
could be computed using~\eqref{eq:D} by a SDP solver.
However, we found it difficult to find its analytical worst-case bound
for the primary sequence $\{\y_i\}$ of OGM,
so~\cite[Prob.~(D$'$)]{kim:17:otc} provided
the following alternate upper bound
on $f(\y_{N+1}) - f(\x_*)$:
\begin{align}
\mathcal{B}_{\mathrm{D'}}(\bmh,N,L,R) :=\;
\min_{\substack{(\bmlam,\bmtau)\in\Lambda, \\ \gamma\in\Reals}}
        \left\{
        \frac{1}{2}LR^2\gamma\;:\;
        \left(\begin{array}{cc}
                \S(\bmh,\bmlam,\bmtau) + \frac{1}{2}\u_N\u_N^\top & \frac{1}{2}\bmtau \\
                \frac{1}{2}\bmtau^\top & \frac{1}{2}\gamma
        \end{array}\right)
        \succeq\Zero
        \right\}
        \label{eq:D_} \tag{D$'$}
,\end{align}
which led to the bound~\eqref{eq:fv_ogm_}
for the primary sequence $\{\y_i\}$ of OGM 
in~\cite{kim:17:otc}.

Similar to~\cite[Lemma~4]{kim:16:ofo}, 
we found a feasible point of~\eqref{eq:D_}
in~\cite[Lemma~3.1]{kim:17:otc},
along with feasible step coefficients \bmh of a~\FO: 
\begin{align}
\hkip &= \begin{cases}
                \frac{1}{t_{i+1}}
                        \paren{2t_k - \sum_{j=k+1}^i\hkj}, & k=0,\ldots,i-1, \\
                1 + \frac{2t_i - 1}{t_{i+1}}, & k=i,
        \end{cases}
       \label{eq:hh_ogm_}
\end{align}
for $t_i$ in~\eqref{eq:t_rule}.
Then,~\cite[Thm.~3.1]{kim:17:otc}
showed the bound~\eqref{eq:fv_ogm_} 
using~\cite[Lemma~3.1]{kim:17:otc}.
The step coefficients~\eqref{eq:hh_ogm} and~\eqref{eq:hh_ogm_}
are identical except the final iteration,
since $t_i$~\eqref{eq:t_rule}
and $\theta_i$~\eqref{eq:theta_rule} are equivalent for $i<N$.
 
We are now done reviewing the portions of papers~\cite{kim:16:ofo,kim:17:otc}
that are the ingredients for 
specifying a parameterized family of algorithms that generalizes OGM
in the next two sections.

\subsection{Feasible points of~\eqref{eq:D} and~\eqref{eq:D_}
for the generalized OGM}

This section
specifies feasible points of~\eqref{eq:D} and~\eqref{eq:D_}
that lead to a generalized version of OGM.
Specifically, the following lemma 
presents additional feasible points 
of~\eqref{eq:D};
this lemma reduces to~\cite[Lemma~4]{kim:16:ofo}
(and the step coefficients~\eqref{eq:hh_ogm} of OGM)
when $\theta_i^2 = \Theta_i$ for all $i$.

\begin{lemma}
\label{lem:gogm}
For the following step coefficients:
\begin{align}
\hkip &= \begin{cases}
                \frac{\theta_{i+1}}{\Theta_{i+1}}
                        \paren{2\theta_k - \sum_{j=k+1}^i\hkj}, 
			& i=0,\ldots,N-1,\;k=0,\ldots,i-1, \\
                1 + \frac{(2\theta_i - 1)\theta_{i+1}}{\Theta_{i+1}}, 
			& i=0,\ldots,N-1,\;k=i,
        \end{cases} 
	\label{eq:hh_gen_ogm}
\end{align}
the choice of variables:
\begin{align}
\gamma = \frac{1}{2}\tau_0,
\quad
\lambda_i &= \Theta_{i-1} \tau_0,\;\; i=1,\ldots,N,
\quad
\tau_i = \begin{cases}
                \frac{2}{\Theta_N}, & i=0, \\
		\theta_i\tau_0 & i=1,\ldots,N-1, \\
                \frac{\theta_N}{2}\tau_0, & i=N,
        \end{cases}
\label{eq:par_gen_ogm}
\end{align}
is a feasible point of~\eqref{eq:D}
for any choice of $\theta_i$ such that
\begin{align}
\theta_0 = 1,
\quad
\theta_i > 0, 
\quad\text{and}\quad
\theta_i^2 \le \Theta_i := \begin{cases}
\sum_{l=0}^i \theta_l, & i=0,\ldots,N-1, \\
2\sum_{l=0}^{N-1}\theta_l + \theta_N, & i=N.
\end{cases}
\label{eq:gen_ogm_rule}
\end{align}
\begin{proof}
See Appendix~\ref{appen2}.
\end{proof}
\end{lemma}

The following lemma also specifies some feasible points of~\eqref{eq:D_};
this lemma reduces to~\cite[Lemma~3.1]{kim:17:otc}
(and~\eqref{eq:hh_ogm_})
when $t_i^2 = T_i$ for all $i$.

\begin{lemma}
\label{lem:gogm_}
For the following step coefficients:
\begin{align}
\hkip &= \begin{cases}
                \frac{t_{i+1}}{T_{i+1}}
                        \paren{2t_k - \sum_{j=k+1}^i\hkj},
                        & i=0,\ldots,N-1,\;k=0,\ldots,i-1, \\
                1 + \frac{(2t_i - 1)t_{i+1}}{T_{i+1}},
                        & i=0,\ldots,N-1,\;k=i,
        \end{cases}
        \label{eq:hh_gen_ogm_}
\end{align}
the choice of variables:
\begin{align}
\gamma = \frac{1}{2}\tau_0,
\quad
\lambda_i &= T_{i-1} \tau_0,\;\; i=1,\ldots,N,
\quad
\tau_i = \begin{cases}
                \frac{1}{T_N}, & i=0, \\
                t_i\tau_0 & i=1,\ldots,N,
        \end{cases}
\label{eq:par_gen_ogm_}
\end{align}
is a feasible point of~\eqref{eq:D_}
for any choice of $t_i$ such that
\begin{align}
t_0 = 1, 
\quad
t_i > 0,
\quad\text{and}\quad
t_i^2 \le T_i := \sum_{l=0}^i t_l.
\label{eq:gen_ogm_rule_}
\end{align}
\begin{proof}
See Appendix~\ref{appen3}.
\end{proof}
\end{lemma}

Similar to the relationship between
the step coefficients~\eqref{eq:hh_ogm} and~\eqref{eq:hh_ogm_},
the step coefficients~\eqref{eq:hh_gen_ogm} and~\eqref{eq:hh_gen_ogm_} 
are identical (when $\theta_i = t_i$ for $i<N$)
except for the final iteration,
implying that the iterates $\{\x_i\}_{i=0}^{N-1}$ 
of the two~{\FO}s
with~\eqref{eq:hh_gen_ogm} and~\eqref{eq:hh_gen_ogm_} are equivalent;
only the final iterate $\x_N$ is different.
The step coefficients~\eqref{eq:hh_ogm} and~\eqref{eq:hh_ogm_} 
lead to computationally efficient equivalent OGM forms;
similarly the next section provides 
computationally efficient generalized forms of OGM
that each correspond to a~\FO 
with either~\eqref{eq:hh_gen_ogm} or~\eqref{eq:hh_gen_ogm_},
and we analyze their cost function worst-case bounds.

\subsection{Generalized OGM}

This section proposes a generalized OGM
using lemmas~\ref{lem:gogm} and~\ref{lem:gogm_}.
The \FO with the step coefficients~\eqref{eq:hh_gen_ogm}
has the following two equivalent efficient generalized forms of OGM,
named GOGM1 and GOGM2,
that reduce to the standard OGM when $\theta_i^2 = \Theta_i$ for all $i$.

\begin{center}
\fbox{
\footnotesize 
\begin{minipage}[t]{0.45\textwidth}
\vspace{-10pt}
\begin{flalign*}
&\text{\bf Algorithm GOGM1} & \\
&\text{Input: } f\in \cF,\; \x_0=\y_0\in\Reals^d, & \\
&	\hspace{30pt} \theta_0=\Theta_0 = 1. & \\
&\text{For } i = 0,\ldots,N-1 & \\
&\quad \y_{i+1} = \x_i - \frac{1}{L}\nabla f(\x_i) & \\
&\quad \text{Choose } \theta_{i+1}>0 & \\
&\quad \text{s.t. }
                \theta_{i+1}^2
                \le \Theta_{i+1}, \text{ where} & \\
&\quad \Theta_{i+1} =
                \begin{cases}
                        \sum_{l=0}^{i+1} \theta_l, & i \le N-2  \\
                        2\sum_{l=0}^{N-1}\theta_l + \theta_N, & i = N-1
                \end{cases} & \\
&\quad \x_{i+1} = \y_{i+1}
                + \frac{(\Theta_i - \theta_i)\theta_{i+1}}{\theta_i\Theta_{i+1}}
			(\y_{i+1} - \y_i) & \\
&\hspace{55pt}    
                + \frac{(2\theta_i^2 - \Theta_i)\theta_{i+1}}{\theta_i\Theta_{i+1}}(\y_{i+1} - \x_i) &
\end{flalign*}
\end{minipage}\vline\hspace{8pt}%
\begin{minipage}[t]{0.45\textwidth}
\vspace{-10pt}
\begin{flalign*}
&\text{\bf Algorithm GOGM2} & \\
&\text{Input: } f\in \cF,\; \x_0=\y_0\in\Reals^d, & \\
&	\hspace{30pt} \theta_0=\Theta_0 = 1. & \\
&\text{For } i = 0,\ldots,N-1 & \\
&\quad \y_{i+1} = \x_i - \frac{1}{L}\nabla f(\x_i) & \\
&\quad \z_{i+1} = \x_0 - \frac{1}{L}\sum_{k=0}^i2\theta_k \nabla f(\x_k) & \\
&\quad \text{Choose } \theta_{i+1}>0 & \\
&\quad \text{s.t. }
                \theta_{i+1}^2
                \le \Theta_{i+1}, \text{ where} & \\
&\quad \Theta_{i+1} =
                \begin{cases}
                        \sum_{l=0}^{i+1} \theta_l, & i \le N-2  \\
                        2\sum_{l=0}^{N-1}\theta_l + \theta_N, & i = N-1
                \end{cases} & \\
&\quad \x_{i+1} = \paren{1 - \frac{\theta_{i+1}}{\Theta_{i+1}}}\y_{i+1}
                + \frac{\theta_{i+1}}{\Theta_{i+1}}\z_{i+1} &
\end{flalign*}
\end{minipage}
} \vspace{5pt}
\end{center}

\begin{proposition}
\label{prop:gogm}
The sequence $\{\x_0,\cdots,\x_N\}$ generated 
by the~\FO with
\eqref{eq:hh_gen_ogm}
is identical to the corresponding sequence
generated by GOGM1 and GOGM2.
\begin{proof}
See Appendix~\ref{appen4}. 
Note that this proof is independent of the choice of
$\theta_i$ and $\Theta_i$.
\end{proof}
\end{proposition}

Because the proof of Prop.~\ref{prop:gogm}
for the~\FO with step coefficients~\eqref{eq:hh_gen_ogm}
is independent of the choice of $\theta_i$ and $\Theta_i$,
it is straightforward to show that
the~\FO with step coefficients~\eqref{eq:hh_gen_ogm_} 
has the following two efficient equivalent forms,
named GOGM1$'$ and GOGM2$'$,
that reduce to~\cite[Alg. OGM1$'$ and OGM2$'$]{kim:17:otc}
when $t_i^2 = T_i$ for all $i$.

\begin{center}
\fbox{
\footnotesize 
\begin{minipage}[t]{0.45\textwidth}
\vspace{-10pt}
\begin{flalign*}
&\text{\bf Algorithm GOGM1$'$} & \\
&\text{Input: } f\in \cF,\; \x_0=\y_0\in\Reals^d, & \\
&	\hspace{30pt} t_0=T_0 = 1. & \\
&\text{For } i = 0,\ldots,N-1 & \\
&\quad \y_{i+1} = \x_i - \frac{1}{L}\nabla f(\x_i) & \\
&\quad \text{Choose } t_{i+1}>0 & \\
&\quad \text{s.t. }
                t_{i+1}^2
                \le T_{i+1} = \sum_{l=0}^{i+1} t_l & \\
&\quad \x_{i+1} = \y_{i+1}
                + \frac{(T_i - t_i)t_{i+1}}{t_iT_{i+1}}
                        (\y_{i+1} - \y_i) & \\
&\hspace{55pt}    
                + \frac{(2t_i^2 - T_i)t_{i+1}}{t_iT_{i+1}}(\y_{i+1} - \x_i)  &
\end{flalign*}
\end{minipage}\vline\hspace{8pt}%
\begin{minipage}[t]{0.45\textwidth}
\vspace{-10pt}
\begin{flalign*}
&\text{\bf Algorithm GOGM2$'$} & \\
&\text{Input: } f\in \cF,\; \x_0=\y_0\in\Reals^d, & \\
&	\hspace{30pt} t_0=T_0 = 1. & \\
&\text{For } i = 0,\ldots,N-1 & \\
&\quad \y_{i+1} = \x_i - \frac{1}{L}\nabla f(\x_i) & \\
&\quad \z_{i+1} = \x_0 - \frac{1}{L}\sum_{k=0}^i2t_k \nabla f(\x_k) & \\
&\quad \text{Choose } t_{i+1}>0 \\
&\quad \text{s.t. }
                t_{i+1}^2
                \le T_{i+1} = \sum_{l=0}^{i+1} t_l & \\
&\quad \x_{i+1} = \paren{1 - \frac{t_{i+1}}{T_{i+1}}}\y_{i+1}
                + \frac{t_{i+1}}{T_{i+1}}\z_{i+1} &
\end{flalign*}
\end{minipage}
} \vspace{5pt}
\end{center}

\noindent
Clearly when $\theta_i = t_i$ for $i<N$,
the primary iterates $\{\y_i\}_{i=0}^N$ 
and the intermediate secondary iterates $\{\x_i\}_{i=0}^{N-1}$ 
of GOGM and GOGM$'$ are equivalent.
Although illustrating two similar algorithms GOGM and GOGM$'$ 
might seem redundant,
presenting both formulations with lemmas~\ref{lem:gogm} and~\ref{lem:gogm_}
completes the story of generalized OGM
here and in Sec.~\ref{sec:pep,sgnorm}.

Using lemmas~\ref{lem:gogm} and~\ref{lem:gogm_},
the following theorem
bounds the cost function decrease of the GOGM and GOGM$'$ iterates.

\begin{theorem}
\label{thm:gogm}
Let $f\;:\;\Reals^d\rightarrow\Reals$ be $\cF$
and let $\y_0,\cdots,\y_N,\x_N\in\Reals^d$ be
generated by GOGM1 and GOGM2. Then for any $1\le i \le N$,
\begin{align}
f(\y_i) - f(\x_*) &\le \frac{LR^2}{4\Theta_{i-1}},
\label{eq:gogm_conv_} \\
f(\x_N) - f(\x_*) &\le \frac{LR^2}{2\Theta_N}
\label{eq:gogm_conv}
.\end{align}
The iterates $\y_0,\cdots,\y_N\in\Reals^d$ 
generated by GOGM1$'$ and GOGM2$'$
also satisfy the bound~\eqref{eq:gogm_conv_}
when $\theta_i = t_i$ for $i<N$.
\begin{proof}
Using Lemma~\ref{lem:gogm_}, 
the~\FO with \bmh in~\eqref{eq:hh_gen_ogm_} 
of GOGM$'$ satisfies
\begin{align}
f(\y_{N+1}) - f(\x_*) \le \mathcal{B}_{\mathrm{D'}}(\bmh,N,L,R)
        = \frac{1}{2}LR^2\gamma
	= \frac{LR^2}{4T_N}
\label{eq:fv_gogm_}
,\end{align}
where $\y_{N+1} = \x_N - \frac{1}{L}\nabla f(\x_N)$.
Since the coefficients \bmh in~\eqref{eq:hh_gen_ogm_}
are recursive and do not depend on a given $N$, 
we can extend~\eqref{eq:fv_gogm_} for all iterations.
By letting $\theta_i = t_i$ for $i<N$,
the bound~\eqref{eq:fv_gogm_} also satisfies
for the iterates $\{\y_i\}$ of GOGM,
as in~\eqref{eq:gogm_conv_}.

Using Lemma~\ref{lem:gogm}, 
the~\FO with the step \bmh in~\eqref{eq:hh_gen_ogm} 
of GOGM satisfies
\begin{align}
f(\x_N) - f(\x_*) \le \mathcal{B}_{\mathrm{D}}(\bmh,N,L,R) 
	= \frac{1}{2}LR^2\gamma 
	= \frac{LR^2}{2\Theta_N} 
,\end{align}
which is equivalent to~\eqref{eq:gogm_conv}.
\end{proof}
\end{theorem}

GOGM and Thm.~\ref{thm:gogm} 
reduce to OGM and its bounds~\eqref{eq:fv_ogm_} and~\eqref{eq:fv_ogm}, 
when $\theta_i^2 = \Theta_i$ for all $i$.
Similar to general forms of FGM in~\cite{chambolle:15:otc,nesterov:05:smo},
the GOGM family includes the choice
$
\theta_i = \begin{cases}
	\frac{i+a}{a}, & i < N, \\	
	\frac{\sqrt{2}(N+a-1)}{a}, & i= N
\end{cases}
$
for any $a\ge2$, because such parameter $\theta_i$ satisfies
the following conditions for GOGM:
\begin{align}
\Theta_i - \theta_i^2 = \frac{(i+1)(i+2a)}{2a} - \frac{(i+a)^2}{a^2}
        = \frac{(a-2)i^2 + a(2a-3)}{2a^2}
        \ge 0,
\label{eq:Thth}
\end{align}
for $i<N$, and
$
\Theta_N - \theta_N^2 = 2\Theta_{N-1} + \theta_N - \theta_N^2
	\ge 2\theta_{N-1}^2 + \theta_N - \theta_N^2	
	= \theta_N \ge 0.
$
Similarly, the GOGM$'$ family includes the choice $t_i=\frac{i+a}{a}$
for any $a\ge2$, which we denote as OGM-$a$.

\begin{corollary}
\label{cor:ogma}
Let $f\;:\;\Reals^d\rightarrow\Reals$ be $\cF$
and let $\y_0,\cdots,\y_N\in\Reals^d$ be
generated by GOGM$'$ with $t_i = \frac{i+a}{a}$ (OGM-$a$)
for any $a\ge2$. Then for any $1\le i \le N$,
\begin{align}
f(\y_i) - f(\x_*) \le \frac{aLR^2}{2i(i+2a-1)}
\label{eq:gogm_cor}
.\end{align}
\begin{proof}
Thm.~\ref{thm:gogm} implies~\eqref{eq:gogm_cor},
since $T_i = \frac{(i+1)(i+2a)}{2a}$
and the condition $T_i - t_i^2\ge0$ satisfies in~\eqref{eq:Thth} 
for any $a\ge2$.
\end{proof}
\end{corollary}

\section{Relaxation and optimization of the gradient form of PEP}
\label{sec:pep,sgnorm}

This section analyzes
a worst-case bound for the gradient of any GOGM (and GOGM$'$)
using the gradient form of PEP.
We use relaxations on the gradient form of PEP
that are similar but slightly different
from those of PEP for the cost function in the previous section.
Using this relaxed PEP,
we prove that FGM has an $O(1/N^{1.5})$ rate for the worst-case gradient decrease,
and analyze the worst-case gradient bound for the GOGM.
Then, we optimize the step coefficients with
respect to the gradient form of PEP
and propose an algorithm named OGM-OG 
that lies in the GOGM family
and that has the best known analytical worst-case bound
for decreasing the gradient norm among the class~\FO.

\subsection{Relaxation for the gradient form of PEP}
\label{sec:pep,grad}

To analyze a worst-case bound on the gradient 
for a~\FO with a given \bmh,
we consider the following gradient-form version of PEP 
that is similar to~\eqref{eq:PEP}:
\begin{align}
\mathcal{B}_{\mathrm{P''}}(\bmh,N,d,L,R) :=\;
& \max_{\substack{f\in\cF, \\
        \x_0,\cdots,\x_N\in\Reals^d, \\
        \x_*\in X_*(f), \\
	||\x_0 - \x_*|| \le R}}
        \min_{i\in\{0,\ldots,N\}}
	||\nabla f(\x_i)||^2
        \label{eq:PEP__} \tag{P$''$} \\
        \st \; &\x_{i+1} = \x_i - \frac{1}{L}\sum_{k=0}^i \hkip \nabla f(\x_k),
                \quad i=0,\ldots,N-1. \nonumber 
\end{align}
Here, we use the smallest gradient norm squared among all iterates
($\min_{i\in\{0,\ldots,N\}}\allowbreak||\nabla f(\x_i)||^2$)
as a criteria, as considered in~\cite[Sec.~4.3]{taylor:17:ssc}.
We could instead consider the final gradient norm squared
($||\nabla f(\x_N)||^2$) as a criteria,
but our proposed relaxation on~\eqref{eq:PEP__} in this section
for such criteria 
provided only an $O(1/N)$ worst-case bound at best
even for the corresponding optimized step coefficients
(results not shown);
we leave studying the gradient form of the tight PEP as future work.

As in~\cite{taylor:17:ssc},
we replace $\min_{i\in\{0,\ldots,N\}} ||\nabla f(\x_i)||^2$ in~\eqref{eq:PEP__}
by $L^2||\x_0-\x_*||^2\alpha$ with the condition 
$\alpha \le \frac{1}{L^2||\x_0-\x_*||^2}||\nabla f(\x_i)||^2 
= \Tr{\G^\top(\u_i\u_i^\top)\G}$ 
for all $i$.
Then, we relax this reformulated~\eqref{eq:PEP__} 
similar to the relaxation from~\eqref{eq:PEP} to~\eqref{eq:pPEP}
with the additional constraint $\Tr{\G^\top\C_N\G} \le \delta_N$ in~\eqref{eq:relPEP}
as follows:
\begin{align}
\mathcal{B}_{\mathrm{P2''}}(\bmh,N,d,L,R) :=\;
&\max_{\substack{\G\in\Reals^{(N+1)d}, \\ \del\in\Reals^{N+1}, \\ \alpha\in\Reals}}
        L^2R^2\alpha
        \label{eq:pPEP__} \tag{P2$''$} \\
        &\st \; \Tr{\G^\top\A_{i-1,i}(\bmh)\G} \le \delta_{i-1} - \delta_i,
                \quad i=1,\ldots,N, \nonumber \\
	&\quad\;\;\;
                \Tr{\G^\top\C_N\G} \le \delta_N, \nonumber \\
        &\quad\;\;\; \Tr{\G^\top\D_i(\bmh)\G + \bmnu \u_i^\top\G} \le -\delta_i,
                \quad i=0,\ldots,N, \nonumber \\
	&\quad\;\;\; \Tr{\G^\top(-\u_i\u_i^\top)\G} \le - \alpha,
		\quad i=0,\ldots,N. \nonumber
\end{align}

Replacing $\max_{\G,\del} L^2R^2\alpha$ 
by $\min_{\G,\del} \{-\alpha\}$ for convenience,
the Lagrangian of the corresponding constrained minimization problem 
with dual variables $\bmlam\in\Reals_+^N$, $\eta\in\Reals_+$, 
$\bmtau\in\Reals_+^{N+1}$, and $\bmbeta = (\beta_0,\cdots,\beta_N)^\top\in\Reals_+^{N+1}$ 
for the first, second, third, and fourth set of constraint
inequalities of~\eqref{eq:pPEP__} respectively becomes
\begin{align}
\cL''(\G,\del,\bmlam,\eta,\bmtau,\bmbeta;\bmh)
        =& - \alpha + \sum_{i=1}^N\lambda_i(\delta_i - \delta_{i-1})
                - \eta\delta_N 
		+ \sum_{i=0}^N\tau_i\delta_i
		+ \sum_{i=0}^N\beta_i\alpha \\
		&\qquad\qquad\qquad 
			+ \Tr{\G^\top\S''(\bmh,\bmlam,\eta,\bmtau,\bmbeta)\G + \bmnu\bmtau^\top\G}
	\nonumber
,\end{align}
where
\begin{align}
\S''(\bmh,\bmlam,\eta,\bmtau,\bmbeta) 
	:= \S(\bmh,\bmlam,\bmtau)
		+ \frac{1}{2}\eta\u_N\u_N^\top
		- \sum_{i=0}^N\beta_i\u_i\u_i^\top
\label{eq:S__}
.\end{align}
Then similar to~\eqref{eq:D}, we have 
the following dual problem of~\eqref{eq:pPEP__}
that one could use to compute an upper bound of the PEP~\eqref{eq:PEP__}
of the smallest gradient norm squared among all iterates
by a numerical SDP solver:
\begin{align}
\mathcal{B}_{\mathrm{D''}}(\bmh,N,L,R) :=\;
&\min_{\substack{(\bmlam,\eta,\bmtau,\bmbeta)\in\Lambda'', \\ \gamma\in\Reals}}
        \left\{
        \frac{1}{2}L^2R^2\gamma\;:\;
        \left(\begin{array}{cc}
                \S''(\bmh,\bmlam,\eta,\bmtau,\bmbeta) & \frac{1}{2}\bmtau \\
                \frac{1}{2}\bmtau^\top & \frac{1}{2}\gamma
        \end{array}\right)
        \succeq\Zero
        \right\}
        \label{eq:D__} \tag{D$''$}
,\end{align}
where
\begin{align}
\Lambda'' := \left\{(\bmlam,\eta,\bmtau,\bmbeta)\in\Reals_+^{3N+3} 
                \;:\;
                \begin{array}{l}
                \tau_0 = \lambda_1, \;\; \lambda_N + \tau_N = \eta, \;\;
                \sum_{i=0}^N \beta_i = 1, \\
                \lambda_i - \lambda_{i+1} + \tau_i = 0, \; i=1,\ldots,N-1
                \end{array}\right\}
\label{eq:Lam__}
.\end{align}

The next two sections use a valid upper bound~\eqref{eq:D__} of~\eqref{eq:PEP__}
for given step coefficients \bmh,
providing an analytical solution to~\eqref{eq:D__}
for the step coefficients \bmh of FGM and GOGM$'$,
superseding the use of a SDP solver.

\subsection{A worst-case bound for the gradient norm of FGM}
\label{sec:pep,sgnorm,fgm}

FGM is equivalent to a~\FO with the step coefficients
\cite[Prop.~1]{kim:16:ofo}:
\begin{align}
h_{i+1,k}
        &= \begin{cases}
                \frac{1}{t_{i+1}}
                        \left(t_k - \sum_{j=k+1}^i h_{j,k}\right), & k=0,\ldots,i-1, \\
                1 + \frac{t_i - 1}{t_{i+1}}, & k=i,
        \end{cases}
        \label{eq:hh_fgm}
\end{align}
for $t_i$ in~\eqref{eq:t_rule}.
The following lemma provides a feasible point of~\eqref{eq:D__}
associated with the step coefficients~\eqref{eq:hh_fgm} of FGM
to provide a worst-case bound for the gradient of FGM.

\begin{lemma}
\label{lem:fgm}
For the step coefficients~\eqref{eq:hh_fgm},
the following choice of variables:
\begin{align}
\gamma &= \tau_0, \quad
\lambda_i = t_{i-1}^2\tau_0, \quad i=1,\ldots,N, \quad
\tau_i =
\begin{cases}
        \left(\frac{1}{2}\sum_{k=0}^N t_k^2\right)^{-1}, & i = 0, \\
        t_i\tau_0, & i = 1,\ldots,N,
\end{cases}
        \label{eq:par_fgm} \\
\eta &= t_N^2 \tau_0, \quad
\beta_i = \frac{1}{2}t_i^2\tau_0, \quad i=0,\ldots,N,
        \label{eq:gam_fgm}
\end{align}
is a feasible point of~\eqref{eq:D__} for $t_i$ in~\eqref{eq:t_rule}.
\begin{proof}
See Appendix~\ref{appen5}.
\end{proof}
\end{lemma}

Using Lemma~\ref{lem:fgm},
the following theorem bounds the gradient norm of the FGM iterates,
proving for the first time an $O(1/N^{1.5})$ rate of decrease.

\begin{theorem}
\label{thm:fgm}
Let $f\;:\;\Reals^d\rightarrow\Reals$ be $\cF$
and let $\y_0,\cdots,\y_N,\x_0,\cdots,\x_N\in\Reals^d$ be
generated by FGM. Then for any $N\ge1$,
\begin{align}
\min_{i\in\{0,\ldots,N+1\}}||\nabla f(\y_i)|| 
&\le \min_{i\in\{0,\ldots,N\}}||\nabla f(\x_i)||  
	\label{eq:fgm_g} \\
&\le \frac{LR}{\sqrt{\sum_{k=0}^Nt_k^2}}
	\le \frac{2\sqrt{3}LR}{\sqrt{(N+1)(N^2 + 6N + 12)}}
	\nonumber
,\end{align}
where $\y_{N+1} = \x_N - \frac{1}{L}\nabla f(\x_N)$.
\begin{proof}
Lemma~\ref{lem:gmono}
implies the first inequality in~\eqref{eq:fgm_g}.
Using Lemma~\ref{lem:fgm}, 
the~\FO with the step coefficients \bmh~\eqref{eq:hh_fgm} of FGM satisfies
\begin{align}
\min_{i\in\{0,\ldots,N\}}||\nabla f(\x_i)||^2 \le \mathcal{B}_{\mathrm{D''}}(\bmh,N,L,R)
        = \frac{1}{2}L^2R^2\gamma
	= \frac{L^2R^2}{\sum_{k=0}^Nt_k^2}
,\end{align}
which is equivalent to~\eqref{eq:fgm_g} using
$ 
\sum_{k=0}^N t_k^2 \ge \sum_{k=0}^N \frac{(k+2)^2}{4}
	= \frac{(N+1)(2N^2 + 13N + 24)}{24}
.$ 
\end{proof}
\end{theorem}

\subsection{A worst-case bound for the gradient norm of GOGM}

Having established
the gradient bound~\eqref{eq:fgm_g}
for FGM,
this section and the next
seek to improve on it
by studying GOGM.
To bound the gradient decrease of GOGM (and GOGM$'$),
the following lemma illustrates one possible set
of feasible points of~\eqref{eq:D__}.

\begin{lemma}
\label{lem:gogm_feas_}
For the step coefficients~\eqref{eq:hh_gen_ogm_},
the following choice of variables:
\begin{align}
\gamma = \frac{1}{2}\tau_0, \quad\!
\lambda_i &= T_{i-1}\tau_0, \quad\! i=1,\ldots,N, \quad\!
\tau_i = 
\begin{cases}
	\left(\sum_{k=0}^N\left(T_k - t_k^2\right)\right)^{-1}, & i = 0, \\
        t_i\tau_0, & i = 1,\ldots,N,
\end{cases}
	\label{eq:par_ogmsg} \\
\eta &= T_N \tau_0, \quad
\beta_i = \left(T_i - t_i^2\right)\tau_0,\quad i=0,\ldots,N
	\label{eq:gam_ogmsg}
\end{align}
is a feasible point of~\eqref{eq:D__} for any choice of $t_i$ and $T_i$
that satisfies~\eqref{eq:gen_ogm_rule_}
and for which there exists some $i$ such that $t_i^2 < T_i$.
\begin{proof}
See Appendix~\ref{appen6}.
\end{proof}
\end{lemma}

Using Lemma~\ref{lem:gogm_feas_},
the following theorem bounds the worst-case gradient norm 
for the iterates of GOGM and GOGM$'$.

\begin{theorem}
\label{thm:gogm_g}
Let $f\;:\;\Reals^d\rightarrow\Reals$ be $\cF$
and let $\y_0,\cdots,\y_N,\x_0,\cdots,\x_N\in\Reals^d$ be
generated by GOGM$'$. Then for any $N\ge1$,
\begin{align}
\min_{i\in\{0,\ldots,N+1\}} ||\nabla f(\y_i)||
\le \min_{i\in\{0,\ldots,N\}} ||\nabla f(\x_i)||
        \le \frac{LR}{2\sqrt{\sum_{k=0}^N\left(T_k - t_k^2\right)}} 
\label{eq:gogm_g}
,\end{align}
where $\y_{N+1} = \x_N - \frac{1}{L}\nabla f(\x_N)$.
The bound~\eqref{eq:gogm_g} can be generalized
to the intermediate iterates $\{\x_i\}_{i=0}^{N-1}$ and $\{\y_i\}_{i=0}^N$
of both GOGM and GOGM$'$ when $\theta_i = t_i$ (for $i<N$).
\begin{proof}
Lemma~\ref{lem:gmono} 
implies the first inequality in~\eqref{eq:gogm_g}.
Using Lemma~\ref{lem:gogm_feas_}, 
\FO with the step coefficients \bmh~\eqref{eq:hh_gen_ogm_} of GOGM$'$ satisfies
\begin{align*}
\min_{i\in\{0,\ldots,N\}} ||\nabla f(\x_i)||^2
\le \mathcal{B}_{\mathrm{D''}}(\bmh,N,L,R)
= \frac{1}{2}L^2R^2\gamma
= \frac{L^2R^2}{4\sum_{k=0}^N\left(T_k - t_k^2\right)}
,\end{align*}
which implies~\eqref{eq:gogm_g}.
Since the iterates of GOGM$'$ are recursive and do not depend on a given $N$,
the bound~\eqref{eq:gogm_g} 
easily generalizes to the intermediate iterates
of GOGM$'$ (and GOGM when $\theta_i = t_i$).
\end{proof}
\end{theorem}

\subsection{Optimizing step coefficients over the gradient form of PEP}
\label{sec:opt,grad,pep}

In search of a~\FO that decreases the gradient norm the fastest,
we optimize the step coefficients 
in terms of the gradient form of the relaxed~\eqref{eq:D__}
by solving the following problem:
\begin{align}
\hat{\bmh}_{\mathrm{D''}}
        := \argmin{\bmh\in\Reals^{N(N+1)/2}} \mathcal{B}_{\mathrm{D''}}(\bmh,N,L,R)
        \tag{HD$''$}
\label{eq:HD__}
.\end{align}
The problem~\eqref{eq:HD__} is bilinear, similar to~\eqref{eq:HD},
and a convex relaxation technique~\cite[Thm.~3]{drori:14:pof}         
makes this problem solvable using numerical methods.

We solved~\eqref{eq:HD__} for many choices of $N$
using a numerical SDP solver~\cite{cvxi,gb08},
and observed that the following choice of $t_i$:
\begin{align}
t_i &= \begin{cases}
                1, & i = 0, \\
                \frac{1+\sqrt{1 + 4t_{i-1}^2}}{2}, & i = 1, \ldots, \fNh - 1, \\
                \frac{N-i+1}{2}, & i = \fNh,\ldots,N,
        \end{cases}
\label{eq:opt_t}
\end{align}
makes the feasible point in Lemma~\ref{lem:gogm_feas_} optimal
for the problem~\eqref{eq:HD__}.
Based on that numerical evidence,
we conjecture that $\hat{\bmh}_{\mathrm{D''}}$ in~\eqref{eq:HD__}
corresponds to the step coefficients~\eqref{eq:hh_gen_ogm_} with
the parameter $t_i$~\eqref{eq:opt_t}.
The $t_i$ factors in~\eqref{eq:opt_t}
start decreasing after $i=\fNh-1$, 
whereas the usual $t_i$ in~\eqref{eq:t_rule}
and $t_i = \frac{i+a}{a}$ for any $a\ge2$
increase with $i$ indefinitely.

In addition, we found numerically that minimizing
the gradient bound~\eqref{eq:gogm_g} of GOGM$'$, \ie,
solving the following constrained quadratic problem:
\begin{align}
\max_{\{t_i\}} \sum_{k=0}^N\paren{\sum_{l=0}^kt_l - t_k^2}
\quad\st\quad t_i \text{ satisfies~\eqref{eq:gen_ogm_rule_} for all } i
\label{eq:quad}
,\end{align}
is equivalent to solving the problem~\eqref{eq:HD__}.
In other words, the solution of~\eqref{eq:quad} numerically
appears equivalent to~\eqref{eq:opt_t},
the (conjectured) solution of~\eqref{eq:HD__}.
The unconstrained maximizer of the cost function of~\eqref{eq:quad}
is $t_i = \frac{N-i+1}{2}$,
and this term
partially appears in the constrained maximizer~\eqref{eq:opt_t}
for $\fNh \le i\le N$.

We denote the resulting GOGM$'$ with~\eqref{eq:opt_t}
as OGM-OG (OG for optimized over gradient).
The following theorem bounds 
the cost function and gradient norm of the OGM-OG iterates.

\begin{theorem}
Let $f\;:\;\Reals^d\rightarrow\Reals$ be $\cF$
and let $\y_0,\cdots,\y_N,\x_0,\cdots,\x_N\in\Reals^d$ be
generated by OGM-OG. Then, 
\begin{align}
f(\y_{N+1}) - f(\x_*) 
        &\le \frac{2LR^2}{(N+2)^2},
\label{eq:ogmog_fv} \\
\min_{i\in\{0,\ldots,N+1\}} ||\nabla f(\y_i)||
&\le \min_{i\in\{0,\ldots,N\}} ||\nabla f(\x_i)|| 
	\le \frac{\sqrt{6}LR}{N\sqrt{N+1}}
\label{eq:ogmog_g}
,\end{align}
where $y_{N+1} = \x_N - \frac{1}{L}\nabla f(\x_N)$.
\begin{proof}
OGM-OG is an instance of GOGM$'$
and thus Thm.~\ref{thm:gogm} implies that OGM-OG satisfies
\begin{align*}
f(\y_{N+1}) - f(\x_*) \le \frac{LR^2}{4T_N},
\end{align*}
which is equivalent to~\eqref{eq:ogmog_fv},
since
\begin{align*}
T_N
        &= T_m + \sum_{l=m+1}^Nt_l
        = t_m^2 + \frac{(N-m)(N-m+1)}{4} \\
	&\ge \frac{(N+3)^2 + N(N+2)}{16}
        = \frac{2N^2 + 8N + 9}{16}
\end{align*}
for $m=\fNh$,
using $m\ge\frac{N-1}{2}$, $N-m \ge\frac{N}{2}$,
and $t_m\ge\frac{m+2}{2}\ge\frac{N+3}{4}$ in~\eqref{eq:t_rule}.

Thm.~\ref{thm:gogm_g} implies~\eqref{eq:ogmog_g},
using 
the above inequalities for $m$,
the equality $t_i^2 = T_i$ for $i\le m$,
and
\begingroup
\allowdisplaybreaks
\begin{align*}
	&\sum_{k=m+1}^N \left(T_k - t_k^2\right)
        = \sum_{k=m+1}^N \left(t_m^2 + \sum_{l=m+1}^k t_l - t_k^2\right) \\
        =& \left(N - m\right)t_m^2
                + \sum_{k=m+1}^N\left(\sum_{l=m+1}^k\frac{N-l+1}{2}
                        - \left(\frac{N-k+1}{2}\right)^2\right) \\
	=& \left(N - m\right)t_m^2
                + \sum_{k'=1}^{N-m}\left(\sum_{l'=1}^{k'}\frac{N-l'-m+1}{2}
                        - \left(\frac{N-k'-m+1}{2}\right)^2\right) \\
	=& \left(N - m\right)t_m^2 \\
        &\quad + \sum_{k=1}^{N-m}\left(\frac{2(N-m+1)k - k(k+1)}{4} 
                        - \frac{(N-m+1)^2 - 2(N-m+1)k + k^2}{4}\right) \\
	=& \left(N - m\right)t_m^2
                + \sum_{k=1}^{N-m}\left(
			- \frac{k^2}{2}
			+ (N-m+3/4)k
			- \frac{(N-m+1)^2}{4}
                        \right) \\
	=& \left(N - m\right)t_m^2
                - \frac{(N-m)(N-m+1/2)(N-m+1)}{6} \\
	&\quad + \frac{(N-m)(N-m+3/4)(N-m+1)}{2} 
		- \frac{(N-m)(N-m+1)^2}{4} \\
        \ge& \frac{(N-m)(m + 2)^2}{4}
		+ \frac{(N-m)^2(N-m+1)}{3}
                - \frac{(N-m)(N-m+1)^2}{4} \\
	\ge& \frac{(N-m)^2(N-m+1)}{3}
        \ge \frac{1}{24}N^2(N+1)
.\end{align*}
\endgroup
\end{proof}
\end{theorem}

\noindent
The gradient bound~\eqref{eq:ogmog_g} of OGM-OG
is asymptotically $\sqrt{2}$-times smaller than that of FGM
in Thm.~\ref{thm:fgm}
and $1.5$-times smaller than that of OGM-$m\!=\!\fNt$
in Thm.~\ref{thm:ogmh}.
Regarding the cost function decrease,
the bound~\eqref{eq:ogmog_fv} of OGM-OG
is asymptotically the same as the bound~\eqref{eq:fv_fgm} of FGM,
and both are twice larger than 
the bounds~\eqref{eq:fv_ogm_} and~\eqref{eq:fv_ogm} of OGM.

\subsection{Decreasing the gradient norm with rate $O(1/N^{1.5})$ 
using GOGM without selecting $N$ in advance}

Although OGM-OG satisfies a small worst-case gradient bound 
with a rate $O(1/N^{1.5})$,
OGM-OG (and FGM-\m and OGM-\m) must select $N$ in advance,
unlike FGM. 
Using Thm.~\ref{thm:gogm_g},
the following corollary shows that OGM-$a$ with $a>2$ 
can decrease the gradient with a rate $O(1/N^{1.5})$
without selecting $N$ in advance.
(Cor.~\ref{cor:ogma} showed that OGM-$a$ algorithm with $a\ge2$
can decrease the cost function with an optimal rate $O(1/N^2)$.)

\begin{corollary}
Let $f\;:\;\Reals^d\rightarrow\Reals$ be $\cF$
and let $\y_0,\cdots,\y_N,\x_0,\cdots,\x_N\in\Reals^d$ be
generated by GOGM$'$ with $t_i = \frac{i+a}{a}$ (OGM-$a$) for any $a\ge2$.
Then for $N\ge1$,
\begin{align}
\min_{i\in\{0,\ldots,N+1\}} ||\nabla f(\y_i)||
&\le \min_{i\in\{0,\ldots,N\}} ||\nabla f(\x_i)|| \label{eq:ogma_g} \\
	&\le \frac{a\sqrt{6}LR}{2\sqrt{N(N+1)\paren{(a-2)N + (3a^2 - 4a - 2)}}}
	\nonumber
,\end{align}
where $\y_{N+1} = \x_N - \frac{1}{L}\nabla f(\x_N)$.
\begin{proof}
Using $T_i = \frac{(i+1)(i+2a)}{2a}$ and~\eqref{eq:Thth},
Thm.~\ref{thm:gogm_g} implies~\eqref{eq:ogma_g} since
\begin{align*}
\sum_{k=0}^N(T_k - t_k^2)
&= \sum_{k=0}^N\frac{(a-2)k^2 + a(2a-3)k}{2a^2} \\
&= \frac{N(N+1)\paren{(a-2)N + (3a^2 - 4a - 2)}}{6a^2}
.\end{align*}
\end{proof}
\end{corollary}

\noindent
OGM-$a$ for any $a>2$ has a gradient bound~\eqref{eq:ogma_g}
that is about $\frac{a}{2\sqrt{a-2}}$-times larger 
than the bound~\eqref{eq:ogmog_g} of OGM-OG.
This constant factor minimizes to $\sqrt{2}$ when $a=4$,
and this OGM-$a\!=\!4$ has a worst-case gradient bound 
that is asymptotically equivalent to the bound~\eqref{eq:fgm_g} of FGM.
Therefore, when one does not want to select $N$ in advance,
both FGM and OGM-$a\!=\!4$ 
(and OGM-$a$ for any $a>2$)
will be useful
for decreasing the gradient with a rate $O(1/N^{1.5})$.

\section{Discussion}
\label{sec:disc}

This section summarizes analytical worst-case bounds of~\FO
discussed in the previous sections.
This section also reports tight numerical worst-case bounds
for exact comparison of algorithms
because many of the analytical bounds are not guaranteed to be tight.

\subsection{Summary of analytical worst-case bounds
on the cost function and gradient norm}

Table~\ref{tab:rate} summarizes
the asymptotic rate
of analytical worst-case bounds of all algorithms 
described in this paper.
As discussed, OGM and OGM-OG have the best known worst-case bounds
for the cost function and gradient decrease respectively in Table~\ref{tab:rate}.
However, since OGM has a slow worst-case rate for the gradient decrease,
other algorithms such as FGM, OGM-\m, OGM-OG, and OGM-$a$
that satisfy both the optimal rate $O(1/N^2)$ for the function decrease
and a fast rate $O(1/N^{1.5})$ for the gradient decrease
could be preferable over OGM when one is interested 
in both the gradient decrease
as well as the function decrease,
particularly when solving dual problems.
In addition, when one does not want to choose $N$ in advance,
FGM and OGM-$a$ could be preferable.

\begin{table}[!h]
\centering
\begin{tabular}{|l|c|c|c|}
\hline
\multirow{2}{*}{Algorithm}
        & \multicolumn{2}{c|}{Asymptotic worst-case bound}
        & Require selecting \\
\cline{2-3}
& Cost function & Gradient norm & $N$ in advance \\
\hline
GM
  & $\frac{1}{4} N^{-1}$ & $\sqrt{2} N^{-1}$ & No  \\ \hline
FGM
  & $2 N^{-2}$  & $2\sqrt{3} N^{-1.5}$ & No  \\ \hline
{\bf OGM}
  & $N^{-2}$  & $\sqrt{2} N^{-1}$ & No  \\ \hline
OGM-$m\!=\!\fNt$  
  & $\frac{9}{4} N^{-2}$ 
	& $\frac{3\sqrt{6}}{2} N^{-1.5}$ & Yes \\ \hline
{\bf OGM-OG}
  & $2 N^{-2}$  & $\sqrt{6} N^{-1.5}$ & Yes \\ \hline
OGM-$a$ ($a > 2$)
  & $\frac{a}{2} N^{-2}$ & $\frac{a\sqrt{6}}{2\sqrt{a-2}} N^{-1.5}$
  & \multirow{2}{*}{No} \\
OGM-$a\!=\!4$
  & $2 N^{-2}$ & $2\sqrt{3} N^{-1.5}$ & \\
\hline
\end{tabular}
\caption{
Asymptotic worst-case bounds
on the cost function $\frac{1}{LR^2}(f(\x_N) - f(\x_*))$
and the gradient norm $\min_{i\in\{0,\ldots,N\}}\frac{1}{LR}||\nabla f(\x_i)||$
of GM, FGM, OGM, OGM-\m, OGM-OG, and OGM-$a$.
(The worst-case cost function bound for OGM-\m in the table corresponds
to the bound for OGM after $\m$ iterations,
because we do not have
an analytical bound for the final iterate.}
\label{tab:rate}
\end{table}
\vspace{-10pt}

\subsection{Tight worst-case bounds
on the cost function and the gradient norm}

Since many worst-case bounds presented in Table~\ref{tab:rate} 
are not guaranteed to be tight, 
we used the code 
in Taylor \etal~\cite{taylor:17:ssc}
(with SDP solvers~\cite{Lofberg2004,sturm:99:us1})
to compare tight (numerical) worst-case bounds 
for $N=1,2,4,10,20,30,40,47,50$.
These numerical worst-case bounds
are guaranteed to be tight, 
\ie, equivalent to the bounds of either~\eqref{eq:PEP} or~\eqref{eq:PEP__},
when the large-scale condition $d\ge N+2$
is satisfied~\cite[Thm.~5]{taylor:17:ssc},
and we assume this condition hereafter.
Tables~\ref{tab:costbound} and~\ref{tab:sgbound} 
provide tight worst-case bounds
for the decrease of the cost function $f(\x_N) - f(\x_*)$
and the gradient norm decrease 
$\min_{i\in\{0,\ldots,N\}}||\nabla f(\x_i)||$ respectively.
Although most of the bounds in Table~\ref{tab:rate} 
are not guaranteed to be tight, 
the worst-case rate formulas in Table~\ref{tab:rate}
are similar to the tight numerical results 
in Tables~\ref{tab:costbound} and~\ref{tab:sgbound},
except that the gradient bounds of OGM-\m in Table~\ref{tab:rate} 
are relatively looser
than those of OGM-\m in Table~\ref{tab:sgbound}. 
In particular,
the tight numerical gradient bound of OGM-\m
is smaller than that of FGM in Table~\ref{tab:sgbound},
which was not expected 
from their known (possibly loose) analytical bounds in Table~\ref{tab:rate}.

\begin{table}[!h]
\centering
\begin{tabular}{|c|l|l|l|l|l|l|}
\hline
$N$ & GM      & FGM      & OGM       & OGM-\m 
					        & OGM-OG   & OGM-$a$ 
										\\
\hline
1   & $6.0$   & $6.0$    & $8.0$     & $6.0$    & $7.3$    & $6.5$    \\
2   & $10.0$  & $11.1$   & $16.2$    & $12.0$   & $13.2$   & $15.1$   \\
4   & $18.0$  & $24.7$   & $39.1$    & $24.2$   & $28.6$   & $32.3$   \\
10  & $42.0$  & $90.7$   & $159.1$   & $86.6$   & $99.9$   & $106.4$  \\
20  & $82.0$  & $283.6$  & $525.1$   & $275.3$  & $310.4$  & $308.9$  \\
30  & $122.0$ & $578.6$  & $1095.6$  & $565.1$  & $604.9$  & $610.9$  \\
40  & $162.0$ & $975.1$  & $1869.2$  & $899.0$  & $1009.9$ & $1012.8$ \\
47  & $190.0$ & $1312.9$ & $2531.1$  & $1227.9$ & $1352.8$ & $1353.6$ \\
50  & $202.0$ & $1472.8$ & $\bm{2845.1}$  & $1374.4$ & $1516.0$ & $1514.6$ \\
\hline
Empi. $O(\cdot)$
    & $N^{-1.0}$ 
    & $N^{-1.9}$ 
    & $N^{-1.9}$
    & $N^{-1.8}$
    & $N^{-1.8}$ 
    & $N^{-1.8}$ \\
\hline
Known $O(\cdot)$   
    & $N^{-1}$   
    & $N^{-2}$  
    & $N^{-2}$ 
    & $N^{-2}$ 
    & $N^{-2}$ 
    & $N^{-2}$  \\
\hline
\end{tabular}
\caption{Tight worst-case bounds on
$\frac{LR^2}{f(\x_N) - f(\x_*)}$,
the reciprocal of the cost function,
of GM, FGM, OGM, OGM-$m\!=\!\fNt$, OGM-OG, and OGM-$a\!=\!4$.
We computed empirical rates
by assuming that the bounds follow the form $bN^{-c}$ with constants $b$ and $c$,
and then by estimating $c$ 
from points $N=47,50$.
Note that the corresponding empirical rates are underestimated
due to its simple modeling of bounds.
}
\label{tab:costbound}
\end{table}

\begin{table}[!h]
\centering
\begin{tabular}{|c|l|l|l|l|l|l|}
\hline
$N$ & GM      & FGM     & OGM     & OGM-\m 
					    & OGM-OG 
						      & OGM-$a$ 
									\\
\hline
1   & $2.0$   & $2.0$   & $2.0$   & $2.0$   & $2.3$   & $1.8$   \\
2   & $3.0$   & $3.3$   & $2.8$   & $3.5$   & $3.7$   & $3.3$   \\
4   & $5.0$   & $5.9$   & $4.4$   & $6.4$   & $6.8$   & $5.7$   \\
10  & $11.0$  & $13.8$  & $8.9$   & $18.0$  & $18.9$  & $15.3$  \\
20  & $21.0$  & $32.8$  & $16.2$  & $43.1$  & $45.4$  & $35.2$  \\
30  & $31.0$  & $56.4$  & $23.4$  & $74.4$  & $78.6$  & $59.4$  \\
40  & $41.0$  & $83.6$  & $30.6$  & $110.7$ & $116.9$ & $87.1$  \\
47  & $48.0$  & $104.7$ & $35.6$  & $138.9$ & $146.6$ & $108.4$ \\
50  & $51.0$  & $114.2$ & $37.7$  & $151.4$ & $\bm{160.0}$ & $118.0$ \\
\hline
Empi. $O(\cdot)$
    & $N^{-1.0}$ & $N^{-1.4}$ & $N^{-0.9}$
    & $N^{-1.4}$ & $N^{-1.4}$ & $N^{-1.4}$ \\
\hline
Known $O(\cdot)$
    & $N^{-1}$ & $N^{-1.5}$ & $N^{-1}$
    & $N^{-1.5}$ & $N^{-1.5}$ & $N^{-1.5}$ \\
\hline
\end{tabular}
\caption{Tight worst-case bounds on
$\frac{LR}{\min_{i\in\{0,\ldots,N\}}||\nabla f(\x_i)||}$,
the reciprocal of the gradient norm,
of GM, FGM, OGM, OGM-$m\!=\!\fNt$, OGM-OG, and OGM-$a\!=\!4$.
Empirical rates were computed
as described in Table~\ref{tab:costbound}.
}
\label{tab:sgbound}
\end{table}

\subsection{Tight worst-case bounds
on the gradient norm at the final iterate}

To be clear, 
Sec.~\ref{sec:pep,sgnorm} and Tables~\ref{tab:rate},~\ref{tab:sgbound}
have focused on analyzing the \emph{smallest} gradient norm among all iterates
using the gradient form of PEP,
whereas the gradient analysis in Sec.~\ref{sec:grad}
considers the \emph{final} gradient 
in addition to the \emph{smallest} gradient among all iterates.
As mentioned before,
we have not yet found 
a relaxation on the \emph{final} gradient form of the PEP
that provides as comparable results
as for the relaxation 
on the \emph{smallest} gradient form of the PEP~\eqref{eq:PEP__}
in Sec.~\ref{sec:pep,sgnorm}.
To complete comparisons on the worst-case gradient bounds, 
Table~\ref{tab:lgbound} uses the code provided by Taylor~\etal~\cite{taylor:17:ssc}
(with SDP solvers~\cite{Lofberg2004,sturm:99:us1})
to compare tight (numerical) worst-case bounds on the \emph{final} gradient
of the~{\FO}s presented in this paper.\footnote{
Table~\ref{tab:lgbound} reports tight worst-case gradient bounds
for both the \emph{final} primary iterate $\y_N$ 
and the \emph{final} secondary iterate $\x_N$ (if necessary),
unlike Tables~\ref{tab:costbound} and~\ref{tab:sgbound}.
We observed that numerical tight worst-case cost function bounds 
on both final iterates $\y_N$ and $\x_N$ 
have similar values
for the algorithms in Table~\ref{tab:costbound} (unlike Table~\ref{tab:lgbound}),
so we did not report the bounds on $\y_N$ for simplicity.
We also did not report numerical tight \emph{smallest} worst-case gradient norm bounds
$\min_{i\in\{0,\ldots,N\}}||\nabla f(\y_i)||$
of the primary iterates $\{\y_i\}$
because
the code provided by Taylor~\etal~\cite{taylor:17:ssc}
does not support computing their values,
unlike that of the secondary iterates $\{\x_i\}$ in Table~\ref{tab:sgbound}.
}
The worst-case \emph{smallest} gradient norm bounds
\eqref{eq:g_gm} and~\eqref{eq:g_ogmh}
of GM and OGM-\m respectively (among algorithms considered)
extend to the \emph{final} gradient bounds.

In Table~\ref{tab:lgbound}, FGM and OGM-$a\!=\!4$ 
have slow $O(1/N)$ tight worst-case bounds on the final gradient,
unlike OGM-$m\!=\!\fNt$ and OGM-OG roughly having $O(1/N^{1.5})$ bounds
for both the smallest and final gradients.
Thm.~\ref{thm:ogmh} has shown that the final gradient of OGM-\m 
satisfies a worst-case rate $O(1/N^{1.5})$,
but this is unknown yet for OGM-OG, which we leave as future work.
We also leave as future work the challenge of developing a~\FO
that has $O(1/N^{1.5})$ or even faster worst-case rates 
for the final gradient decrease
that are lower than those of OGM-\m and OGM-OG,
possibly without requiring to choose $N$ in advance.

\begin{table}[!h]
\centering
\begin{tabular}{|c|l|l|l|l|l|l|l|l|l|}
\hline
\multirow{2}{*}{$N$}  
   & \multirow{2}{*}{GM}      
	     & \multicolumn{2}{c|}{FGM}     
			       & \multicolumn{2}{c|}{OGM}     
				                & OGM  
							  & OGM 
						       		    & \multicolumn{2}{c|}{OGM-$a$} \\
\cline{3-6}\cline{9-10}
   &        & $\y_N$
                     & $\x_N$ & $\y_N$ 
                                       & $\x_N$ & -\m     
							  & -OG     & $\y_N$ 
                                                                             & $\x_N$ \\
\hline
1  & $2.0$  & $2.0$  & $2.0$  & $2.0$  & $2.0$  & $2.0$   & $2.3$   & $2.0$  & $1.8$  \\
2  & $3.0$  & $3.0$  & $3.3$  & $3.2$  & $2.8$  & $3.5$   & $3.7$   & $3.6$  & $3.3$  \\
4  & $5.0$  & $5.8$  & $5.9$  & $5.5$  & $4.4$  & $6.4$   & $6.8$   & $6.8$  & $5.1$  \\
10 & $11.0$ & $15.1$ & $8.2$  & $11.9$ & $8.9$  & $18.0$  & $18.9$  & $15.9$ & $8.7$  \\
20 & $21.0$ & $25.1$ & $13.1$ & $22.2$ & $16.2$ & $43.1$  & $44.4$  & $26.3$ & $13.8$ \\
30 & $31.0$ & $35.1$ & $18.2$ & $32.4$ & $23.4$ & $74.4$  & $74.1$  & $36.3$ & $18.8$ \\
40 & $41.0$ & $45.2$ & $23.2$ & $42.5$ & $30.6$ & $110.7$ & $107.0$ & $46.3$ & $23.8$ \\
47 & $48.0$ & $52.2$ & $26.7$ & $49.6$ & $35.6$ & $138.9$ & $131.6$ & $53.3$ & $27.3$ \\
50 & $51.0$ & $55.3$ & $28.2$ & $52.6$ & $37.7$ & $\bm{151.4}$ & $142.6$ & $56.3$ & $28.8$ \\
\hline
Empi. $O(\cdot)$ 
   & $N^{-1.0}$ & \multicolumn{2}{c|}{$N^{-0.9}$} & \multicolumn{2}{c|}{$N^{-0.9}$}
   & $N^{-1.4}$ & $N^{-1.3}$ & \multicolumn{2}{c|}{$N^{-0.9}$} \\
\hline
Known $O(\cdot)$
   & $N^{-1}$ & \multicolumn{2}{c|}{$N^{-1}$} & \multicolumn{2}{c|}{$N^{-1}$} 
   & $N^{-1.5}$ & $N^{-1}$ & \multicolumn{2}{c|}{$N^{-1}$} \\ 
\hline
\end{tabular}
\caption{
Tight worst-case bounds on 
$\frac{LR}{||\nabla f(\x_N)||}$
$\paren{\text{and } \frac{LR}{||\nabla f(\y_N)||}}$,
the reciprocal of the final gradient norm,
of GM, FGM, OGM, OGM-$m\!=\!\fNt$, OGM-OG, and OGM-$a\!=\!4$.
Empirical rates were computed
as described in Table~\ref{tab:costbound}.
The known bounds of OGM-OG and OGM-$a\!=\!4$ 
are derived based on Sec.~\ref{sec:grad},
where the empirical bounds of OGM-OG
are comparable to the bounds of OGM-$m\!=\!\fNt$ 
with known rate $O(1/N^{1.5})$.
}
\label{tab:lgbound}
\end{table}

\subsection{Non-optimality
of OGM-OG in terms of the worst-case gradient bound}

Because OGM is optimal in terms of the function decrease
when $d \ge N+1$~\cite{drori:17:tei},
one might hope that the OGM-OG would achieve 
the optimal worst-case bound
in terms of the gradient decrease,
since OGM-OG is also derived by optimizing the step coefficients
over the gradient form of relaxed PEP.
However, the OGM-OG is apparently not optimal as explained next.

Taylor~\etal~\cite{taylor:17:ssc}
numerically studied an optimal fixed-step GM
using their tight PEP
in terms of both the cost function and gradient decrease.
In other words, 
they searched for an optimal step $h$ of GM:
\begin{align*}
\x_{i+1} = \x_i - \frac{h}{L}\nabla f(\x_i)
\end{align*}
for $i=0,\ldots,N-1$ and a given $N$
with respect to either $f(\x_N) - f(x_*)$ or $||f(\x_N)||$. 
In the special case of $N=1$,
Taylor~\etal~\cite{taylor:17:ssc} numerically conjectured that
the step size $h=1.5$ is optimal in terms of the cost function decrease.
The corresponding GM is equivalent to OGM for $N=1$,
and this (numerically) confirms the optimality of OGM
\cite{drori:17:tei}
for $N=1$.
They also numerically conjectured that
the optimal step size of GM for $N=1$ 
in terms of the gradient decrease is $h=\sqrt{2}$
with a worst-case bound  
\begin{align}
||\nabla f(\x_1)|| \le \frac{LR}{\sqrt{2}+1}
\approx \frac{LR}{2.4}
.\end{align}
However, OGM-OG for $N=1$ reduces to GM with $h=\frac{4}{3}\approx1.3$
with a bound $\frac{LR}{2.3}$ in Tables~\ref{tab:sgbound} and~\ref{tab:lgbound},
implying that OGM-OG is not optimal even for $N=1$
based on the numerical evidence in~\cite{taylor:17:ssc}.

This analysis for $N=1$ illustrates that there is still room for improvement
in accelerating the worst-case rate of first-order methods 
in terms of gradients,
which we leave as future work possibly with a tighter relaxation 
on the gradient form of PEP.
In addition, 
we leave as future work studying
the optimal worst-case bound for the gradient decrease of first-order methods
building upon~\cite{drori:17:tei,nemirovsky:92:ibc}, 
and developing a~\FO that achieves such optimal bound.
Nevertheless, 
the OGM-OG is the best known~\FO for decreasing the gradient norm
among the class~\FO,
and will be useful when decreasing the gradient is key.

\comment{
\subsection{
Generalizing OGM for nonsmooth composite convex problems}

One can extend some aspects
of the approaches for generalizing OGM described in this paper
to other optimization algorithms and problems.
One direction we have already taken in~\cite{kim:18:ala}
aims to improve the fast iterative shrinkage/thresholding algorithm (FISTA)
\cite{beck:09:afi}
(that reduces to FGM for smooth convex problems)
for nonsmooth composite convex problems.
Naturally, this paper and~\cite{kim:18:ala}
use some similar approaches,
but they are different in the following two aspects.
First, the methods in~\cite{kim:18:ala}
when simplified to the smooth case correspond to a generalization of FGM
that differs from the methods proposed here.
Second,~\cite{kim:18:ala} uses a relaxation for PEP
that is looser than the relaxation in this paper
for the smooth convex case,
and thus the worst-case analysis in~\cite{kim:18:ala}
does not lead to an $O(1/N^{1.5})$ worst-case gradient rate for FGM here.
Nevertheless, the FISTA generalizations in~\cite{kim:18:ala}
provide new insights on first-order methods 
for nonsmooth composite convex problems,
as this paper does for the smooth convex problem.
Finding a straightforward extension of this paper
to nonsmooth composite convex problems
remains as future work.
}

\section{Conclusion}
\label{sec:conc}

We generalized the formulation of OGM 
and analyzed its worst-case bounds on the function value and gradient,
using the cost function form and the gradient form of relaxed PEP. 
We then proposed OGM-OG
by optimizing the step coefficients of~\FO using a relaxed PEP
with respect to the gradient,
similar to the development of the (optimal) OGM.
To the best of our knowledge,
the worst-case bound on the gradient of the OGM-OG
is the best known analytical worst-case bound
for decreasing the gradient norm among the class~\FO.

However, this OGM-OG is not optimal
for decreasing the gradient norm,
and further accelerating the worst-case rate of~\FO in terms of the gradient
possibly with a tight relaxation on the gradient form of PEP
is a possible research direction.
On the other hand,
deriving an optimal worst-case bound
for the gradient norm of first-order methods,
similar to that for the function decrease~\cite{drori:17:tei}
will be useful.
Nonetheless, the proposed OGM-OG (and OGM-$a$)
may be useful when one finds minimizing gradients important,
particularly in dual problems. 
In addition, we used the proposed gradient form of PEP
to show that FGM decreases the (smallest) gradient
with a rate $O(1/N^{1.5})$,
implying that FGM is comparable in a big-O sense to OGM-\m, OGM-OG and OGM-$a$
for the gradient decrease.

Our analysis considers unconstrained smooth convex minimization;
extending such gradient norm worst-case analysis
to constrained problems or nonsmooth composite convex problems
is a natural direction to pursue,
which is studied for FGM (or FISTA~\cite{beck:09:afi})
by the authors~\cite{kim:18:ala}.
In addition,
extending the analyses on the general form of FGM
in~\cite{attouch:18:fco,chambolle:15:otc,su:16:ade}
to GOGM
is a possible research direction.
Lastly,
investigating a new relaxation of the PEP approach
that allows adaptive step size 
such as backtracking line-search
or exact line-search~\cite{deklerk:17:otw,drori:18:efo}
is of interest.

\section*{Software}
\url{https://gitlab.eecs.umich.edu/michigan-fast-optimization}
has
Matlab codes for the algorithms considered and
the SDP approaches
in Sec.~\ref{sec:opt,grad,pep} and Sec.~\ref{sec:disc}.

\appendix

\section{Proof of Thm.~\ref{thm:g_ogm_tight}}
\label{appen1}

Due to~\eqref{eq:g_ogm},
we have
\begin{align*}
\max_{\substack{f\in\cF, \\ \x_* \in X_∗(f), \\ ||\x_0 - \x_*||\le R}} \min_{i\in\{0,\ldots,N\}}
                ||\nabla f(\x_i)||
        \le \max_{\substack{f\in\cF, \\ \x_* \in X_∗(f), \\ ||\x_0 - \x_*||\le R}} ||\nabla f(\x_N)||
        \le \frac{LR}{\theta_N}
,\end{align*}
and the rest of the proof shows that the following inequality holds
\begin{align}
\frac{LR}{\theta_N}
&= \min_{i\in\{0,\ldots,N\}} ||\nabla \phi(\x_i)||
= ||\nabla \phi(\x_N)|| \label{eq:g_ogm_low} \\
&\le
\max_{\substack{f\in\cF, \\ \x_* \in X_∗(f), \\ ||\x_0 - \x_*||\le R}} \min_{i\in\{0,\ldots,N\}}
                ||\nabla f(\x_i)||
                \le \max_{\substack{f\in\cF, \\ \x_* \in X_∗(f), \\ ||\x_0 - \x_*||\le R}} ||\nabla f(\x_N)||
	\nonumber
,\end{align}
which then implies~\eqref{eq:g_ogm_tight}
with $\theta_N \ge \frac{N+1}{\sqrt{2}}$~\eqref{eq:theta_rule}.

Starting from $\x_0 = R\bmnu$, where \bmnu is a unit vector,
we first use induction to show that
the following iterates:
\begin{align}
\x_i = (-1)^i\frac{1}{\theta_i}R\bmnu, \quad i=0,\ldots,N,
\label{eq:qpoint}
\end{align}
correspond to the iterates of OGM applied to $\phi(\x)$.
We use~\cite[Prop.~4]{kim:16:ofo} that the sequence generated by OGM
is identical to the sequence generated by~\FO with
\begin{align}
h_{i+1,k} = \begin{cases}
        \frac{\theta_i-1}{\theta_{i+1}}h_{i,k}, & k=0,\ldots,i-2, \\
        \frac{\theta_i-1}{\theta_{i+1}}(h_{i,i-1} - 1), & k=i-1, \\
        1 + \frac{2\theta_i-1}{\theta_{i+1}}, & k=i,
\end{cases}
\label{eq:hh_ogm1}
\end{align}
for $i=0,\ldots,N-1$.

Assuming that~\eqref{eq:qpoint} holds for $i < N$,
we have
\begingroup
\allowdisplaybreaks
\begin{align*}
&\x_{i+1} = \x_i - \frac{1}{L}\sum_{k=0}^i\hkip\nabla\phi(\x_k) \\
        =& \x_i 
	- \frac{1}{L}\paren{1 + \frac{2\theta_i-1}{\theta_{i+1}}}\nabla\phi(\x_i)
                - \frac{1}{L}\sum_{k=0}^{i-1}\frac{\theta_i - 1}{\theta_{i+1}}
                        h_{i,k}\nabla\phi(\x_k)
                + \frac{1}{L}\frac{\theta_i-1}{\theta_{i+1}}\nabla\phi(\x_{i-1}) \\
        =& \frac{1 - 2\theta_i}{\theta_{i+1}}\x_i
                + \frac{\theta_i-1}{\theta_{i+1}}(\x_i - \x_{i-1})
                + \frac{\theta_i-1}{\theta_{i+1}}\x_{i-1} \\
        =& - \frac{\theta_i}{\theta_{i+1}}\x_i
        = (-1)^{i+1}\frac{1}{\theta_{i+1}}R\bmnu,
\end{align*}
\endgroup
using~\eqref{eq:hh_ogm1} and $\nabla\phi(\x) = L\x$.
Therefore, after $N$ iterations of OGM we have:
\begin{align}
\min_{i\in\{0,\ldots,N\}} ||\nabla \phi(\x_i)||
= ||\nabla\phi(\x_N)||
= \bigg|\bigg|\nabla\phi\paren{(-1)^N\frac{1}{\theta_N}R\bmnu}\bigg|\bigg|
= \frac{LR}{\theta_N}
\label{eq:g_ogm_iter}
,\end{align}
which is equivalent to~\eqref{eq:g_ogm_low}.
The first equality of~\eqref{eq:g_ogm_iter}
holds since OGM monotonically decreases the gradient norm of $\phi(\x)$,
\ie $||\nabla\phi(\x_i)|| = \frac{LR}{\theta_i}$
monotonically decreases as $i$ increases.

\section{Proof of Lemma~\ref{lem:gogm}}
\label{appen2}

It is obvious that $(\bmlam, \bmtau)$ in~\eqref{eq:par_gen_ogm}
is in $\Lambda$~\eqref{eq:Lam}.
Using~\eqref{eq:ABCDF} and~\eqref{eq:S},
we have
\begingroup
\allowdisplaybreaks
\begin{align}
\S(\bmh,\bmlam,\bmtau)
&= \begin{cases}
\frac{1}{2}\paren{(\lambda_i + \tau_i)\hki + \tau_i\sum_{j=k+1}^{i-1}\hkj},
        & i=2,\ldots,N,\;k=0,\ldots,i-2, \\
\frac{1}{2}\paren{(\lambda_i + \tau_i)\hki - \lambda_i},
        & i=1,\ldots,N,\;k=i-1, \\
\lambda_{i+1}, & i=0,\ldots,N-1,\;k=i, \\
\frac{1}{2}, & i=N,\;k=i,
\end{cases}
\label{eq:SS}
\end{align}
and inserting~\eqref{eq:hh_gen_ogm} and~\eqref{eq:par_gen_ogm},
yields
\begin{align*}
&\S(\bmh,\bmlam,\bmtau) \\
=& \begin{cases}
\frac{1}{2}\paren{\Theta_i\tau_0\frac{\theta_i}{\Theta_i}
        \paren{2\theta_k - \sum_{j=k+1}^{i-1}\hkj}
        + \theta_i\tau_0\sum_{j=k+1}^{i-1}\hkj},
	& i=2,\ldots,N-1,\;k=0,\ldots,i-2, \\
\frac{1}{2}\paren{\frac{\Theta_N}{2}\tau_0\frac{\theta_N}{\Theta_N}
        \paren{2\theta_k - \sum_{j=k+1}^{N-1}\hkj}
        + \frac{\theta_N}{2}\tau_0\sum_{j=k+1}^{N-1}\hkj},
        & i=N,\;k=0,\ldots,i-2, \\
\frac{1}{2}\paren{\Theta_i\tau_0
        \paren{1 + \frac{(2\theta_{i-1} - 1)\theta_i}{\Theta_i}}
        - \Theta_{i-1}\tau_0},
        & i=1,\ldots,N-1,\;k=i-1, \\
\frac{1}{2}\paren{\frac{\Theta_N}{2}\tau_0
        \paren{1 + \frac{(2\theta_{N-1} - 1)\theta_N}{\Theta_N}}
        - \Theta_{N-1}\tau_0},
        & i=N,\;k=i-1, \\
\Theta_i\tau_0,
        & i=0,\ldots,N-1,\;k=i, \\
\frac{\Theta_N}{4}\tau_0,
        & i=N,\;k=i,
\end{cases} \\
=& \begin{cases}
\theta_i\theta_k\tau_0,
        & i=1,\ldots,N-1,\;k=0,\ldots,i-1, \\
\frac{\theta_N\theta_k}{2}\tau_0,
        & i=N,\;k=0,\dots,i-1, \\
\Theta_i\tau_0,
        & i=1,\ldots,N-1,\;k=i, \\
\frac{\Theta_N}{4}\tau_0,
        & i=N,\;k=i,
\end{cases}
\end{align*}
\endgroup
for $\theta_i$ and $\Theta_i$ in~\eqref{eq:gen_ogm_rule}.
Then, using~\eqref{eq:par_gen_ogm} and~\eqref{eq:gen_ogm_rule},
we show the feasibility condition of~\eqref{eq:D}:
\begin{align*}
\left(\begin{array}{cc}
                \S(\bmh,\bmlam,\bmtau) & \frac{1}{2}\bmtau \\
                \frac{1}{2}\bmtau^\top & \frac{1}{2}\gamma
        \end{array}\right)
        = \paren{\diag{\vTheta - \vtheta^2} + \vtheta\vtheta^\top}\tau_0
        \succeq\Zero
,\end{align*}
where
$\vtheta = \paren{\theta_0,\cdots,\theta_{N-1},\frac{\theta_N}{2},\frac{1}{2}}^\top$
and
$\vTheta = \paren{\Theta_0,\cdots,\Theta_{N-1},\frac{\Theta_N}{4},\frac{1}{4}}^\top$.

\section{Proof of Lemma~\ref{lem:gogm_}}
\label{appen3}

It is obvious that $(\bmlam, \bmtau)$ in~\eqref{eq:par_gen_ogm_}
is in $\Lambda$~\eqref{eq:Lam}.
\begingroup
\allowdisplaybreaks
Inserting~\eqref{eq:hh_gen_ogm_} and~\eqref{eq:par_gen_ogm_}
to~\eqref{eq:SS}
yields
\begin{align*}
&\S(\bmh,\bmlam,\bmtau) + \frac{1}{2}\u_N\u_N^\top \\
=& \begin{cases}
\frac{1}{2}\paren{T_i\tau_0\frac{t_i}{T_i}
        \paren{2t_k - \sum_{j=k+1}^{i-1}\hkj}
        + t_i\tau_0\sum_{j=k+1}^{i-1}\hkj},
        & i=2,\ldots,N-1,\;k=0,\ldots,i-2, \\
\frac{1}{2}\paren{T_i\tau_0
        \paren{1 + \frac{(2t_{i-1} - 1)t_i}{T_i}}
        - T_{i-1}\tau_0},
        & i=1,\ldots,N,\;k=i-1, \\
T_i\tau_0,
        & i=0,\ldots,N,\;k=i
\end{cases} \\
=& \begin{cases}
t_it_k\tau_0,
        & i=1,\ldots,N,\;k=0,\ldots,i-1, \\
T_i\tau_0,
        & i=1,\ldots,N,\;k=i, \\
\end{cases}
\end{align*}
\endgroup
for $t_i$ and $T_i$ in~\eqref{eq:gen_ogm_rule_}.
Then, using~\eqref{eq:par_gen_ogm_} and~\eqref{eq:gen_ogm_rule_},
we show the feasibility condition of~\eqref{eq:D_}:
\begin{align*}
\left(\begin{array}{cc}
                \S'(\bmh,\bmlam,\bmtau) & \frac{1}{2}\bmtau \\
                \frac{1}{2}\bmtau^\top & \frac{1}{2}\gamma
        \end{array}\right)
        = \paren{\diag{\vT - \vt^2} + \vt\vt^\top}\tau_0
        \succeq\Zero
,\end{align*}
where
$\vt = \paren{t_0,\cdots,t_N,\frac{1}{2}}^\top$
and
$\vT = \paren{T_0,\cdots,T_N,\frac{1}{4}}^\top$.

\section{Proof of Prop.~\ref{prop:gogm}}
\label{appen4}

The proof consists of three propositions
and they follow the derivations in
\cite[Prop.~3, 4 and 5]{kim:16:ofo} respectively.
Note that this proof is independent of the choice
of $\theta_i$ and $\Theta_i$.

\begin{proposition}
\label{prop:hh_gogm}
The step coefficient~\eqref{eq:hh_gen_ogm} satisfies
the following recursive relationship
\begin{align}
\hkip = \begin{cases}
                \frac{(\Theta_i-\theta_i)\theta_{i+1}}{\theta_i\Theta_{i+1}}\hki
                        & k=0,\ldots,i-2, \\
                \frac{(\Theta_i-\theta_i)\theta_{i+1}}{\theta_i\Theta_{i+1}}(h_{i,i-1} - 1),
                        & k=i-1, \\
                1 + \frac{(2\theta_i - 1)\theta_{i+1}}{\Theta_{i+1}},
                        & k=i,
        \end{cases}
\label{eq:hh_gen_ogm1}
\end{align}
for $i=0,\ldots,N-1$.
\begin{proof}
We use the notation $\hki'$ for the coefficients~\eqref{eq:hh_gen_ogm}
to distinguish from~\eqref{eq:hh_gen_ogm1}.
It is obvious that $\hiip' = \hiip, i=0,\ldots,N-1$,
and we clearly have
\begin{align*}
\himip' &= \frac{\theta_{i+1}}{\Theta_{i+1}}\paren{2\theta_{i-1} - \himi'}
        = \frac{\theta_{i+1}}{\Theta_{i+1}}
        \paren{2\theta_{i-1}
        - \paren{1 + \frac{(2\theta_{i-1}-1)\theta_i}{\Theta_i}}} \\
        &= \frac{(2\theta_{i-1}-1)(\Theta_i - \theta_i)\theta_{i+1}}
                {\Theta_i\Theta_{i+1}}
        = \frac{(\Theta_i - \theta_i)\theta_{i+1}}{\theta_i\Theta_{i+1}}(\himi - 1)
        = \himip
.\end{align*}
We next use induction by assuming $\hkip' = \hkip$
for $i=0,\ldots,n-1,\;k=0,\ldots,i$. We then have
\begin{align*}
\hknp' &= \frac{\theta_{n+1}}{\Theta_{n+1}}
                \paren{2\theta_k - \sum_{j=k+1}^n\hkj'}
        = \frac{\theta_{n+1}}{\Theta_{n+1}}
                \paren{2\theta_k - \sum_{j=k+1}^{n-1}\hkj' - \hkn'} \\
        &= \frac{\theta_{n+1}}{\Theta_{n+1}}
                \paren{\frac{\Theta_n}{\theta_n}\hkn' - \hkn'}
        = \frac{(\Theta_n - \theta_n)\theta_{n+1}}{\theta_n\Theta_{n+1}}\hkn
        = \hknp
.\end{align*}
\end{proof}
\end{proposition}

\begin{proposition}
\label{prop:gogm1}
The sequence $\{\x_0,\ldots,\x_N\}$ generated by~\FO with~\eqref{eq:hh_gen_ogm1}
is identical to the corresponding sequence generated 
by GOGM1.
\begin{proof}
We use induction, and for clarity, we use the notation
$\x_0',\cdots,\x_N'$ for~\FO with~\eqref{eq:hh_gen_ogm1}.
It is obvious that $\x_0' = \x_0$, and we have
\begin{align*}
\x_1' &= \x_0' - \frac{1}{L}h_{1,0}\nabla f(\x_0')
        = \x_0 - \frac{1}{L}\paren{1 + \frac{(2\theta_0 - 1)\theta_1}{\Theta_1}}\nabla f(\x_0) \\
        &= \y_1 + \frac{(\Theta_0 - \theta_0)\theta_1}{\theta_0\Theta_1}(\y_1 - \y_0)
                + \frac{(2\theta_0^2 - \Theta_0)\theta_1}{\theta_0\Theta_1}(\y_1 - \x_0)
        = \x_1
.\end{align*}
Assuming $\x_i' = \x_i$ for $i=0,\ldots,n$, we then have
\begingroup
\allowdisplaybreaks
\begin{align*}
&\x_{n+1}'
= \x_n' - \frac{1}{L}h_{n+1,n}\nabla f(\x_n')
        - \frac{1}{L}h_{n+1,n-1}\nabla f(\x_{n-1}')
        - \frac{1}{L}\sum_{k=0}^{n-2}h_{n+1,k}\nabla f(\x_k') \\
=& \x_n - \frac{1}{L}\paren{1 + \frac{(2\theta_n-1)\theta_{n+1}}
                {\Theta_{n+1}}}\nabla f(\x_n)
        - \frac{1}{L}\frac{(\Theta_n-\theta_n)\theta_{n+1}}{\theta_n\Theta_{n+1}}
                (h_{n,n-1} - 1)\nabla f(\x_{n-1}) \\
        &\quad - \frac{1}{L}\sum_{k=0}^{n-2}
                \frac{(\Theta_n - \theta_n)\theta_{n+1}}{\theta_n\Theta_{n+1}}
                h_{n,k}\nabla f(\x_k) \\
=& \y_{n+1} - \frac{1}{L}\frac{(2\theta_n^2 - \Theta_n)\theta_{n+1}}
		{\theta_n\Theta_{n+1}}\nabla f(\x_n) \\
	&\quad - \frac{1}{L}\frac{(\Theta_n-\theta_n)\theta_{n+1}}{\theta_n\Theta_{n+1}}
        \paren{\nabla f(\x_n) - \nabla f(\x_{n-1}) + \sum_{k=0}^{n-1}h_{n,k}\nabla f(\x_k)} \\
=& \y_{n+1} + \frac{(2\theta_n^2 - \Theta_n)\theta_{n+1}}
		{\theta_n\Theta_{n+1}}(\y_{n+1} - \x_n) \\
        &\quad + \frac{(\Theta_n-\theta_n)\theta_{n+1}}{\theta_n\Theta_{n+1}}
                \paren{-\frac{1}{L}\nabla f(\x_n) + \frac{1}{L}\nabla f(\x_{n-1})
                + \x_n - \x_{n-1}} \\
=& \y_{n+1} + \frac{(\Theta_n-\theta_n)\theta_{n+1}}{\theta_n\Theta_{n+1}}(\y_{n+1} - \y_n)
        + \frac{(2\theta_n^2 - \Theta_n)\theta_{n+1}}{\theta_n\Theta_{n+1}}(\y_{n+1} - \x_n)
= \x_{n+1}
.\end{align*}
\endgroup
\end{proof}
\end{proposition}

\begin{proposition}
\label{prop:gogm2}
The sequence $\{\x_0,\cdots,\x_N\}$ generated by~\FO with~\eqref{eq:hh_gen_ogm}
is identical to the corresponding sequence generated 
by GOGM2.
\begin{proof}
We use induction, and for clarity, we use the notation
$\x_0',\cdots,\x_N'$ for~\FO with~\eqref{eq:hh_gen_ogm}.
It is obvious that $\x_0' = \x_0$, and we have
\begin{align*}
\x_1' &= \x_0' - \frac{1}{L}h_{1,0}\nabla f(\x_0')
        = \x_0 - \frac{1}{L}\paren{1 + \frac{(2\theta_0 - 1)\theta_1}{\Theta_1}}\nabla f(\x_0) \\
        &= \paren{1 - \frac{\theta_1}{\Theta_1}}\paren{\x_0 - \frac{1}{L}\nabla f(\x_0)}
		+ \frac{\theta_1}{\Theta_1}\paren{\x_0 - \frac{1}{L}\nabla f(\x_0)
		- \frac{1}{L}(2\theta_0 - 1)\nabla f(\x_0)} \\
	&= \paren{1 - \frac{\theta_1}{\Theta_1}}\y_1
		+ \frac{\theta_1}{\Theta_1}\z_1
        = \x_1
.\end{align*}
Assuming $\x_i' = \x_i$ for $i=0,\ldots,n$, we then have
\begingroup
\allowdisplaybreaks
\begin{align*}
&\x_{n+1}'
= \x_n' - \frac{1}{L}h_{n+1,n}\nabla f(\x_n')
        - \frac{1}{L}\sum_{k=0}^{n-1}h_{n+1,k}\nabla f(\x_k') \\
=& \x_n - \frac{1}{L}\paren{1 + \frac{(2\theta_n-1)\theta_{n+1}}
                {\Theta_{n+1}}}\nabla f(\x_n)
	- \frac{1}{L}\sum_{k=0}^{n-1}
                \frac{\theta_{n+1}}{\Theta_{n+1}}
		\paren{2\theta_k - \sum_{j=k+1}^n h_{j,k}}
		\nabla f(\x_k) \\
=& \paren{1 - \frac{\theta_{n+1}}{\Theta_{n+1}}}\paren{\x_n - \frac{1}{L}\nabla f(\x_n)} \\
	&\quad + \frac{\theta_{n+1}}{\Theta_{n+1}}
	\paren{\x_n
		- \frac{1}{L}\sum_{k=0}^n 2\theta_k\nabla f(\x_k) 
		+ \frac{1}{L}\sum_{k=0}^{n-1}\sum_{j=k+1}^n h_{j,k}\nabla f(\x_k)} \\
=& \paren{1 - \frac{\theta_{n+1}}{\Theta_{n+1}}}\paren{\x_n - \frac{1}{L}\nabla f(\x_n)}
        + \frac{\theta_{n+1}}{\Theta_{n+1}}
        \paren{\x_0 - \frac{1}{L}\sum_{k=0}^n 2\theta_k\nabla f(\x_k)} \\
=& \paren{1 - \frac{\theta_{n+1}}{\Theta_{n+1}}}\y_{n+1}
        + \frac{\theta_{n+1}}{\Theta_{n+1}}\z_{n+1}
.\end{align*}
\endgroup
\end{proof}
\end{proposition}

\section{Proof of Lemma~\ref{lem:fgm}}
\label{appen5}

It is obvious that $(\bmlam, \eta, \bmtau, \bmbeta)$
in~\eqref{eq:par_fgm} and~\eqref{eq:gam_fgm}
is in $\Lambda''$~\eqref{eq:Lam__}.
Using~\eqref{eq:S__} and~\eqref{eq:SS}, we have
\begin{align}
&\S''(\bmh,\bmlam,\eta,\bmtau,\bmbeta) \label{eq:SS__} \\
=& \begin{cases}
\frac{1}{2}\paren{(\lambda_i + \tau_i)\hki + \tau_i\sum_{j=k+1}^{i-1}\hkj},
        & i=2,\ldots,N,\;k=0,\ldots,i-1, \\
\frac{1}{2}\paren{(\lambda_i + \tau_i)\hki - \lambda_i},
        & i=1,\ldots,N,\;k=i-1, \\
\lambda_{i+1} - \beta_i, & i=0,\ldots,N-1,\;k=i, \\
\eta - \beta_N, & i=N,\;k=i,
\end{cases}
\nonumber
\end{align}
and inserting~\eqref{eq:hh_fgm},~\eqref{eq:par_fgm}, and~\eqref{eq:gam_fgm} yields
\begin{align*}
&\S''(\bmh,\bmlam,\eta,\bmtau,\bmbeta) \\
=& \begin{cases}
\frac{1}{2}\paren{t_i^2\tau_0\frac{1}{t_i}
        \paren{t_k - \sum_{j=k+1}^{i-1}\hkj}
        + t_i\tau_0\sum_{j=k+1}^{i-1}\hkj},
        & i=2,\ldots,N,\;k=0,\ldots,i-1, \\
\frac{1}{2}\paren{t_i^2\tau_0
        \paren{1 + \frac{t_{i-1} - 1}{t_i}}
        - t_{i-1}^2\tau_0},
        & i=1,\ldots,N,\;k=i-1, \\
\frac{1}{2}t_i^2\tau_0,
        & i=0,\ldots,N,\;k=i,
\end{cases} \\
=& \frac{1}{2}t_it_k\tau_0,\quad i=0,\ldots,N,\;k=0,\ldots,i,
\end{align*}
for $t_i$ in~\eqref{eq:t_rule}.
Then, using~\eqref{eq:par_fgm}, we finally show that
the feasibility condition of~\eqref{eq:D__} holds:
\begin{align*}
\left(\begin{array}{cc}
                \S''(\bmh,\bmlam,\eta,\bmtau,\bmbeta) & \frac{1}{2}\bmtau \\
                \frac{1}{2}\bmtau^\top & \frac{1}{2}\gamma
        \end{array}\right)
        = \frac{1}{2}\vt\vt^\top\tau_0
        \succeq\Zero
,\end{align*}
where
$\vt = \paren{t_0,\cdots,t_N,1}^\top$.

\section{Proof of Lemma~\ref{lem:gogm_feas_}}
\label{appen6}

It is obvious that $(\bmlam, \eta, \bmtau, \bmbeta)$ 
in~\eqref{eq:par_ogmsg} and~\eqref{eq:gam_ogmsg}
is in $\Lambda''$.
Inserting~\eqref{eq:hh_gen_ogm_},~\eqref{eq:par_ogmsg} and~\eqref{eq:gam_ogmsg}
to~\eqref{eq:SS__} yields
\begin{align*}
&\S''(\bmh,\bmlam,\eta,\bmtau,\bmbeta) \\
=& \begin{cases}
\frac{1}{2}\paren{T_i\tau_0\frac{t_i}{T_i}
        \paren{2t_k - \sum_{j=k+1}^{i-1}\hkj}
        + t_i\tau_0\sum_{j=k+1}^{i-1}\hkj},
        & i=2,\ldots,N,\;k=0,\ldots,i-1, \\
\frac{1}{2}\paren{T_i\tau_0
        \paren{1 + \frac{(2t_{i-1} - 1)t_i}{T_i}}
        - T_{i-1}\tau_0},
        & i=1,\ldots,N,\;k=i-1, \\
T_i\tau_0 - \paren{T_i - t_i^2}\tau_0,
        & i=0,\ldots,N,\;k=i,
\end{cases} \\
=& t_it_k\tau_0,\quad i=0,\ldots,N,\;k=0,\ldots,i,
\end{align*}
for $t_i$ and $T_i$ in~\eqref{eq:gen_ogm_rule_}.
Then, using~\eqref{eq:par_ogmsg}, we finally show that
the feasibility condition of~\eqref{eq:D__} holds:
\begin{align*}
\left(\begin{array}{cc}
                \S''(\bmh,\bmlam,\eta,\bmtau,\bmbeta) & \frac{1}{2}\bmtau \\
                \frac{1}{2}\bmtau^\top & \frac{1}{2}\gamma
        \end{array}\right)
        = \vt\vt^\top\tau_0
        \succeq\Zero
,\end{align*}
where
$\vt = \paren{t_0,\cdots,t_N,\frac{1}{2}}^\top$.

\section*{Acknowledgements}
The authors would like to thank the anonymous referees
for very useful comments that improved the quality of this paper.

        \bibliographystyle{siamplain}
\bibliography{master,mastersub}

\end{document}

%% file: def,command.tex
\input def,gen

\input def,math

\usepackage{jf-accent}
\usepackage{jf-fun}
\usepackage{jf-names}

\newcommand{\st} {\xmath{\text{s.t.}\:}}

%% file: def,gen.tex
\providecommand{\etal}		{\emph{et al\@.}\xspace}
\providecommand{\ie}		{\emph{i.e\@.}\xspace}
\providecommand{\eg}		{\emph{e.g\@.}\xspace}

\providecommand{\myurl}[1][]	{\texttt{web.eecs.umich.edu/$\sim$fessler#1}\xspace}

\providecommand{\onweb}[1]	{Available from \myurl.}

\long\def\comment#1{}

\providecommand{\bcent}		{\begin{center}}
\providecommand{\ecent}		{\end{center}}
\providecommand{\benum}		{\begin{enumerate}}
\providecommand{\eenum}		{\end{enumerate}}
\providecommand{\bitem}		{\begin{itemize}}
\providecommand{\eitem}		{\end{itemize}}
\providecommand{\bvers}		{\begin{verse}}
\providecommand{\evers}		{\end{verse}}
\providecommand{\btab}		{\begin{tabbing}}	
\providecommand{\etab}		{\end{tabbing}}

\newcounter{blist}
\providecommand{\blistmark}	{\makebox[0pt]{$\bullet$}}
\providecommand{\blistitemsep}	{0pt}
\providecommand{\blist}[1][]	{%
\begin{list}{\blistmark}{%
\usecounter{blist}%
\setlength{\itemsep}{\blistitemsep}%
\setlength{\parsep}{0pt}%
\setlength{\parskip}{0pt}%
\setlength{\partopsep}{0pt}%
\setlength{\topsep}{0pt}%
\setlength{\leftmargin}{1.2em}%
\setlength{\labelsep}{0.5\leftmargin}
\setlength{\labelwidth}{0em}%
#1}
}
\providecommand{\elist}		{\end{list}}

\providecommand{\blistitemsep}	{0pt}
\providecommand{\bjfenum}[1][]	{%
\begin{list}{\bcolor{\arabic{blist}.} }{%
\usecounter{blist}%
\setlength{\itemsep}{\blistitemsep}%
\setlength{\parsep}{0pt}%
\setlength{\parskip}{0pt}%
\setlength{\partopsep}{0pt}%
\setlength{\topsep}{0pt}%
\setlength{\leftmargin}{0.0em}%
\setlength{\labelsep}{1.0\leftmargin}
\setlength{\labelwidth}{0pt}%
#1}
}

\newcounter{blistAlph}
\providecommand{\blistAlph}[1][]
{\begin{list}{\makebox[0pt][l]{\Alph{blistAlph}.}}{%
\usecounter{blistAlph}%
\setlength{\itemsep}{0pt}\setlength{\parsep}{0pt}%
\setlength{\parskip}{0pt}\setlength{\partopsep}{0pt}%
\setlength{\topsep}{0pt}%
\setlength{\leftmargin}{1.2em}%
\setlength{\labelsep}{1.0\leftmargin}
\setlength{\labelwidth}{0.0\leftmargin}#1}%
}

\newcounter{blistRoman}
\providecommand{\blistRoman}[1][]
{\begin{list}{\Roman{blistRoman}.}{%
\usecounter{blistRoman}%
\setlength{\itemsep}{0.5em}\setlength{\parsep}{0pt}%
\setlength{\parskip}{0pt}\setlength{\partopsep}{0pt}%
\setlength{\topsep}{0pt}%
\setlength{\leftmargin}{4em}%
\setlength{\labelsep}{0.4\leftmargin}
\setlength{\labelwidth}{0.6\leftmargin}#1}%
}

%% file: def,math.tex
\usepackage{bbm} 
\providecommand{\jfbbm}[1]	{\xmath{\mathbbm{#1}}} 

\providecommand{\reals}		{\jfbbm{R}}

\providecommand{\floor}[1]	{\xmath{\left\lfloor #1 \right\rfloor}}

\providecommand{\inprod}[2]	{\xmath{\mathop{\langle #1,\, #2 \rangle}\nolimits}}

\let\equivsave\equiv
\def\equiv{\xmath{\equivsave}}

\providecommand{\ba}[1]		{\left[ \begin{array}{#1}}
\providecommand{\ea}		{\end{array} \right]}
\providecommand{\be}		{\begin{equation}}
\providecommand{\ee}[1]		{\label{#1}\end{equation}}
\providecommand{\bea}		{\begin{eqnarray}}
\providecommand{\eea}[1]	{\label{#1}\end{eqnarray}}
\providecommand{\beas}		{\begin{eqnarray*}}
\providecommand{\eeas}		{\end{eqnarray*}}
\providecommand{\beals}[1][1]	{\begin{alignat*}{#1}}	
\providecommand{\eeals}		{\end{alignat*}}

\providecommand{\berr}[2]{
\bgroup
\renewcommand{\theequation}{#1}
\be
#2
\ee{e,#1}
\egroup
\ignorespaces
}

\providecommand{\bearr}[2]{
\bgroup
\renewcommand{\theequation}{#1}
\bea
#2
\eea{e,#1}
\egroup
\ignorespaces
}

\providecommand{\inmath}	{\ensuremath}
\providecommand{\xmath}[1]	{\inmath{#1}\xspace}
\providecommand{\bmath}[1]	{\xmath{\bm{#1}}}	

\providecommand{\paren}[1]	{\xmath{\left(#1\right)}}

\providecommand{\braces}[1]	{\xmath{\left\{#1\right\}}}

\providecommand{\xfrac}[2]	{\xmath{\frac{#1}{#2}}}

\providecommand{\Frac}[2]	{\xmath{{#1}/{#2}}}
\providecommand{\pfrac}[3][]	{\xmath{\!\paren{#1\frac{#2}{#3}}}}

\providecommand{\wordbrace}[3][]{\,\inmath{\mathrm{#2}#1\!\braces{#3}}}

\providecommand{\diag}[1]	{\wordbrace{diag}{#1}}

\providecommand{\sinc}		{\Ofun{\mathrm{sinc}}}

